\documentclass[11pt]{article}
\usepackage[russian, english]{babel}
\usepackage[cp1251]{inputenc}
\usepackage[russian]{babel}
\usepackage{amssymb,amsmath}
\begin{document}
\begin{center}
\textbf{\Large OPTIMAL QUADRATURE FORMULAS WITH POSITIVE
COEFFICIENTS IN $L_2^{(m)}(0,1)$ SPACE}\\
Kh.M.Shadimetov, A.R.Hayotov\\
\end{center}

\large
\begin{abstract}
In the Sobolev space $L_2^{(m)}(0,1)$ optimal quadrature formulas
with the nodes (1.5) are investigated. For optimal coefficients
explicit form are obtained and norm of the error functional is
calculated. In particular, by choosing parameter $\eta_0$ in (1.5)
the optimal quadrature formulas with positive coefficients are
obtained and compared with well known optimal formulas.
\end{abstract}



\section{Introduction. Statement of the problem}

It is known, that numerical integration formulae, or quadrature
formulae, are methods for the approximate evaluation of definite
integrals. They are needed for the computation of those integrals
for which either the antiderivative of the integrand cannot be
expressed in terms of elementary functions or for which the
integrand is available only at discrete points, for example from
experimental data. In addition and even more important, quadrature
formulae provide a basic and important tool for the numerical
solution of differential and integral equations.

There are various methods in the theory of quadrature, which allow
us approximately calculate integrals with the help of finite
number of values of integrand. The present paper also is devoted
to one of such methods, i.e. to construction of optimal quadrature
formulas for approximate evaluation of definite integrals in the
space $L_2^{(m)}(0,1)$ equipped with the norm
$$
\|\varphi(x)\|_{L_2^{(m)}(0,1)}=
\left\{\int\limits_0^1(\varphi^{(m)}(x))^2dx\right\}^{1/2}.
$$

Consider a quadrature formula
$$
\int\limits_0^1 {\varphi (x)dx \cong \sum\limits_{\beta  = 0}^N
{C_\beta  \varphi (x_\beta  )} }\eqno (1.1)
$$
with the error functional
$$
\ell(x) = \varepsilon _{[0,1]} (x) - \sum\limits_{\beta  =
0}^N {C_\beta  \delta (x - x_\beta  )}\eqno (1.2)
$$
in the Sobolev space $L_2^{(m)} (0,1)$. Here $C_\beta$ and
$x_\beta$ are the coefficients and the nodes of the quadrature
formula (1.1), respectively, $\delta (x)$ is the Dirac's
delta-function, $\varepsilon _{[0,1]}(x)$ is the indicator of the
interval $[0,1]$.

In order that the error functional (1.2) is defined on the space
$L_2^{(m)}(0,1)$ it is necessary performing following conditions
(see [28])
$$
(\ell(x),x^{\alpha})=0,\ \ \alpha=0,1,2,...,m-1.\eqno (1.3)
$$

The difference
$$
(\ell (x),\varphi (x)) = \int\limits_0^1 {\varphi (x)dx -
\sum\limits_{\beta  = 0}^N {C_\beta  \varphi (x_\beta  )} }  =
\int\limits_{ - \infty }^\infty  {\ell (x)\varphi (x)dx}.\eqno
(1.4)
$$
is called by \emph{the error} of the quadrature formula (1.1).

By Cauchy-Schwartz inequality
$$
|(\ell(x),\varphi(x))|\leq \|\varphi(x)|L_2^{(m)}(0,1)\|\cdot
\|\ell(x)|L_2^{(m)*}(0,1)\|
$$
the error (1.4) of the formula (1.1) is estimated with the help of
norm
$$
\left\|{\ell (x)|L_2^{(m)*} } \right\| = \mathop {\sup
}\limits_{\left\| {\varphi (x)|L_2^{(m)} } \right\| = 1} |(\ell
(x),\varphi (x))|.
$$
of the error functional (1.2). Consequently, estimation of the
error (1.4) of the quadrature formula (1.1) on functions of the
space $L_2^{(m)}(0,1)$ is reduced to finding of norm of the error
functional $\ell(x)$ in the conjugate space $L_2^{(m)*}(0,1)$.

Clearly, that norm of the error functional $\ell (x)$ depends on
the coefficients $C_{\beta}$ and the nodes $x_{\beta}$. The
problem of finding of minimum of norm of the error functional
$\ell(x)$ by coefficients $C_{\beta}$ and by nodes $x_{\beta}$, is
called by \emph{S.M.Nikolskii problem}, and obtained formula is
called \emph{optimal quadrature formula in the sense of
Nikolskii}. This problem first considered by S.M.Nikolskii [15],
and continued by many authors, see e.g. [2-4,16,32] and references
therein. Minimization of norm of the error functional $\ell(x)$ by
coefficients $C_{\beta}$ when the nodes are fixed is called
\emph{Sard's problem}. And obtained formula is called
\emph{optimal quadrature formula in the sense of Sard}. First this
problem investigated by A.Sard [17].

There are several methods of construction of optimal quadrature
formulas in the sense of Sard such as spline method, $\varphi-$
function method (see e.g. [23], [2]) and Sobolev's method which
based on construction of discrete analogue of a linear
differential operator (see e.g. [26]). In the different spaces,
based on these methods, the Sard's problem investigated by many
authors, see, for example, [1,2,4,6,7,10-14,19-26,28-31] and
references therein.

Furthermore, in works [11,19,21,31] were obtained explicit
formulas for coefficients of optimal quadrature formulas for any
$m$ and for any number $N$ of the nodes $x_{\beta}$ in the space
$L_2^{(m)}$.

By I.J.Schoenberg and S.D.Silliman [25] investigated the Sard's
problem for the case $N \to \infty $ in the space $L_2^{(m)}$. In
[25] an algorithm for finding of optimal coefficients was given
with the help of spline of degree $2m - 1$. In the cases $m =
2,3,...,7$ the coefficients are calculated using a Computer. There
was especially noted, that when $m=7$ among optimal coefficients
the negative coefficient $B_4^{(7)}$ appears. It appears that by
increasing $m$ the number of negative coefficients increases. This
is confirmed by computations of the optimal coefficients for
$m\leq 30$ in [30].

But, as known that in applications the optimal quadrature formulas
with positive coefficients play very important role.

Naturally, it is arises a question: can we obtain with some way
the optimal quadrature formulas with positive coefficients in the
Sobolev space $L_2^{(m)}(0,1)$?

The main objective of the present paper is, using the Sobolev's
method, to construct of optimal quadrature formulas in the sense
of Sard in the space $L_2^{(m)}(0,1)$ with the nodes
$$
x_i  = \eta _i h,\,\,\,\,x_{N - i}  = 1 - \eta _i h,\,\,\,\,\,i =
\overline {0,t - 1} ,\,\,\,0 \le \eta _0  < \eta _1  < ... < \eta
_{t - 1} ,
$$
$$
x_\beta   = h\beta ,\,\,\,\,t \le \beta  \le N - t,\,\,\,\,\,h =
{1 \over N},\ \ \ t =\left\{
\begin{array}{ll}
{{m \over 2}}& \mbox{ when } m \mbox{ is even},\\
\left[ {{m \over 2}}\right]+1& \mbox{ when } m  \mbox{ is odd},
\end{array}
\right. \eqno (1.5)
$$
where $\left[ a \right]$ is integer part of the number $a$ and
choosing the parameters $\eta_i$, $i=\overline{0,t-1}$ to obtain
the optimal quadrature formulas of the form (1.1) with positive
coefficients. This means to find the coefficients $C_{\beta}$
which satisfy following equality
$$
\left\| {\mathop \ell \limits^ \circ  (x)|L_2^{(m)*} } \right\|  =
\mathop {\inf }\limits_{C_{\beta}} \left\| {\ell (x)|L_2^{(m)*} }
\right\| \eqno     (1.6)
$$
for the nodes (1.5).

Thus, in order to construct of optimal quadrature formula in the
sense of Sard with nodes (1.5) in the space $L_2^{(m)}(0,1)$ we
need consequently to solve following problems.

\textbf{Problem 1.} \emph{Find norm of the error functional  $\ell
(x)$ of quadrature formulas of the form (1.1) in the space
$L_2^{(m)*} (0,1)$.}

\textbf{Problem 2.} \emph{Find the coefficients $C_{\beta}$ which
satisfy the equality (1.6) with nodes (1.5).}

It is known [26,28,29],  that square of the norm of the error
functional (1.2) for arbitrary fixed $x_\beta$ has following form
$$
\left\| {\ell(x) |L_2^{(m)*} (0,1)} \right\|^2  = ( - 1)^m \left(
{\sum\limits_{\beta  = 0}^N {\sum\limits_{\gamma  = 0}^N {C_\beta
C_\gamma  {{|x_\beta   - x_\gamma  |^{2m - 1} } \over {2 \cdot (2m
- 1)!}} - } } } \right.
$$
$$
 - \left. {2\sum\limits_{\beta  = 0}^N {C_\beta  {{x_\beta ^{2m}  +
 (1 - x_\beta ^{} )^{2m} } \over {2 \cdot (2m)!}}}  + {1 \over {(2m + 1)!}}} \right).
\eqno (1.7)
$$

Obviously, that the square of the norm (1.7) of the error
functional $\ell(x)$ is multidimensional function with respect to
the coefficients $C_{\beta}$. Moreover, the error functional
$\ell(x)$ satisfies the conditions (1.3). In  [28,29] taking into
account these facts the Lagrange function with conditions (1.3) is
constructed for finding \emph{condition minimum} of the (1.7) and
differentiating that function by $C_\beta$ and by $\lambda
_\alpha$ (where $\lambda _\alpha$ are Lagrange factors) following
system of linear equations was obtained
$$
\sum\limits_{\gamma  = 0}^N {C_\gamma  {{|x_\beta   - x_\gamma
|^{2m - 1} } \over {2 \cdot (2m - 1)!}} + \sum\limits_{\alpha  =
0}^{m - 1} {\lambda _\alpha  x_\beta ^\alpha   = {{x_\beta ^{2m} +
(1 - x_\beta  )^{2m} } \over {2\left( {2m} \right)!}},\ \ \ \
x_\beta \in [0,1],} } \eqno (1.8)
$$
$$
\sum\limits_{\gamma  = 0}^N {C_\gamma  x_\gamma ^\alpha   = {1
\over {\alpha  + 1}},\,\,\,\,\,\alpha  = \overline {0,m - 1}.}
\eqno (1.9)
$$

It was proved in [28,29] that this system has unique solution and
this solution gives minimum to the expression (1.7). This means
that square of norm of the error functional $\ell(x)$ being
quadratic function of the coefficients $C_{\beta}$ has unique
minimum in concrete value of
$C_{\beta}=\stackrel{\circ}{C}_{\beta}$. The quadrature formula
with the coefficients $\stackrel{\circ}{C}_{\beta}$ when the nodes
are fixed is called \emph{optimal quadrature formula in the sense
of Sard} and the coefficients $\stackrel{\circ}{C}_{\beta}$ are
called \emph{optimal one}.

Below for convenience the optimal coefficients
$\stackrel{\circ}{C}_{\beta}$ we shall remain as $C_{\beta}$.

Thus, the problem 1 is already solved by S.L.Sobolev in the space
$L_2^{(m)}$ and the problem 2 is reduced to the system of linear
equations for optimal coefficients. Should be noted, that in
[28,29] the problem 1 solved for multidimensional case, i.e. for
cubature formulas.

In the present paper we will solve the system (1.8)-(1.9), i.e. we
will solve the problem 2.

\section{Definitions and known formulas}

In this section we give some definitions and formulas which are
necessary in the proof of the main results.

The \emph{Euler-Frobenius polynomials} $E_k (x)$ , $k = 1,2,...$
are defined by following formula (see, e.g. [26])
$$
E_k (x) = \frac{{(1 - x)^{k + 2} }}{x}\left( {x\frac{d}{{dx}}}
\right)^k \frac{x}{{(1 - x)^2 }},\eqno     (2.1)
$$
$E_0 (x) = 1$.

For Euler-Frobenius polynomials following identity holds
$$
E_k (x) = x^k E_k \left( {\frac{1}{x}} \right),  \eqno (2.2)
$$
and also following theorem is take placed

\textbf{Theorem 2.1} [20]. \emph{Polynomial $P_k (x)$ which is
determined by formula}
    $$
P_k (x) = (x - 1)^{k + 1} \sum\limits_{i = 0}^{k + 1}
{\frac{{\Delta ^i 0^{k + 1} }}{{(x - 1)^i }}}\eqno (2.3)
$$
\emph{is the Euler-Frobenius polynomial (2.1) of degree $k$, i.e.
$P_k (x) = E_k(x)$.}

Following formula is valid [9]:
$$
\sum\limits_{\gamma  = 0}^{n - 1} {q^\gamma  \gamma ^k  =
\frac{1}{{1 - q}}\sum\limits_{i = 0}^k {\left( {\frac{q}{{1 - q}}}
\right)^i \Delta ^i 0^k  - \frac{{q^n }}{{1 - q}}\sum\limits_{i =
0}^k {\left( {\frac{q}{{1 - q}}} \right)^i \Delta ^i \gamma ^k
|_{\gamma  = n} ,} } }\eqno    (2.4)
$$
where $\Delta^i\gamma^k$ is finite difference of order $i$ of
$\gamma^k$, $\Delta^i0^k=\Delta^i\gamma^k|_{\gamma=0}$.\\
At last we give following well known formulas from [8]
$$
\sum\limits_{\gamma  = 0}^{\beta  - 1} {\gamma ^k  =
\sum\limits_{j = 1}^{k + 1} {\frac{{k!\,B_{k + 1 - j} }}{{j!\,(k +
1 - j)!}}\,\beta ^j ,} }\eqno    (2.5)
$$
where $B_{k + 1 - j} $ are Bernoulli numbers,
    $$
\Delta ^\alpha  x^\nu   = \sum\limits_{p = 0}^\nu  {\left(
{\begin{array}{c}
   \nu   \\
   p  \\
\end{array}} \right)\Delta ^\alpha  } 0^p x^{\nu  - p}.
\eqno (2.6)
$$

\section{Auxiliary results}

\subsection{Lemmas}

In the proofs of the main results we need following lemmas.

\textbf{Lemma 3.1.} [20] {\it The Optimal coefficients of
quadrature formulas of  the form (1.1) with the error functional
(1.2) and the nodes (1.5) in space $L_2^{(m)}(0,1)$ have following
form
$$
C_\beta   = h\left( {1 + \sum\limits_{k = 1}^{m - 1} {d_k \left(
{q_k^\beta   + q_k^{N - \beta } } \right)} } \right) \eqno (3.1)
$$
when $\beta  = t,t + 1,...,N - t.$ Here $d_k$ are unknown
parameters, $q_k$ are roots of the Euler-Frobenius polynomial
$E_{2m - 2}(q)$, which $\left| {q_k } \right| < 1$.}

\textbf{Lemma 3.2.} \textit{Following identity holds
    $$
\sum\limits_{i = 1}^\alpha  {{{ - q_k^{N - s + 1}  + ( - 1)^i
q_k^{s + i} } \over {\left( {q_k  - 1} \right)^{i + 1} }}} \Delta
^i 0^\alpha   = ( - 1)^\alpha  \sum\limits_{i = 1}^\alpha
{{{q_k^{s + 1}  + ( - 1)^{i + 1} q_k^{N - s + i} } \over {\left(
{q_k  - 1} \right)^{i + 1} }}} \Delta ^i 0^\alpha,\eqno (3.2)
$$
where $\alpha  = 1,2,...,m - 1,\,\,\,\,s = 0,1,...,$
 $q_k$ are roots of the Euler-Frobenius polynomial $E_{2m - 2}(q)$.}

\textbf{Lemma 3.3.} \textit{The coefficients of optimal quadrature
formulas of the form (1.1) with the error functional (1.2) and the
nodes (1.5) in the space $L_2^{(m)}(0,1)$ satisfy following system
$$
\sum\limits_{\beta  = 0}^{t - 1} {C_\beta  \eta _\beta ^\alpha   =
h\left( {\sum\limits_{\beta  = 1}^{t - 1} {\beta ^\alpha   +
{{0^\alpha  } \over 2} + \sum\limits_{k = 1}^{m - 1} {d_k
\sum\limits_{i = 0}^\alpha  {{{( - 1)^i q_k^{t + i}  - q_k^{N - t
+ 1} } \over {(q_k  - 1)^{i + 1} }}\Delta ^i t^\alpha  } } } }
\right)}\eqno (3.3)
$$
here $ \alpha  = \left\{ \begin{array}{l}
  0,2,4,...,m - 2\,\,\,\,\,when\,\,\,\,m - \,even, \hfill \cr
  0,2,4,...,m - 1\,\,\,\,\,when\,\,\,\,m - \,odd, \hfill \cr
  \end{array}  \right.
$
 $0^\alpha   = \left\{
 \begin{array}{l}
  1,\,\,\,\,\alpha  = 0, \hfill \cr
  0,\,\,\,\alpha  \ne 0, \hfill \cr
  \end{array}  \right.
$\\
$t =\left\{
\begin{array}{ll}
 {{m \over 2}}& \mbox{ when } m \mbox{ is even},\\
\left[ {{m \over 2}}\right]+1& \mbox{ when } m  \mbox{ is odd},
\end{array}
\right.$ where $\left[ a \right] $ is integer part of the number
$a$, $d_k$  are unknown parameters, $q_k$ are roots of the
Euler-Frobenius polynomial $E_{2m - 2}(q)$, $\left| {q_k } \right|
< 1 $.}

\textbf{Lemma 3.4.} \emph{Let $m > n$ and $m,n \in N$ then
following formula for binomial coefficients is true}
$$
C_m^{m - 1} C_{m - 1}^n  - C_m^{m - 2} C_{m - 2}^n  + ... + ( -
1)^{m - n - 2} C_m^{n + 1} C_{n + 1}^n  + ( - 1)^{m - n - 1} C_m^n
C_n^n  = C_m^n .\eqno (3.4)
$$

\textbf{Lemma 3.5.} [28] \textit{Following identity is true
$$
y^k  = \sum\limits_{i = 1}^k {{{y^{\left[ k \right]} } \over
{i!}}\Delta ^i 0^k ,}\eqno (3.5)
$$
where $y^{\left[ k \right]}  = y(y - 1)(y - 2)...(y - k + 1) $,
$\Delta ^i 0^k   = \sum\limits_{l = 1}^i {( - 1)^{i - l} C_i^l l^k
}$.}

\subsection{Proofs of lemmas}

Lemma 3.1 is proved in [20].

\textbf{Proof of lemma 3.2.} We will denote the left and the right
hand sides of the identity (3.2) by $A_1$ and $A_2$, respectively,
i.e.
$$
A_1  = \sum\limits_{i = 1}^\alpha  {{{ - q_k^{N - s + 1}  + ( -
1)^i q_k^{s + i} } \over {\left( {q_k  - 1} \right)^{i + 1} }}}
\Delta ^i 0^\alpha  ,\,\,\,\,\,\,A_2  = ( - 1)^\alpha
\sum\limits_{i = 1}^\alpha  {{{q_k^{s + 1}  + ( - 1)^{i + 1}
q_k^{N - s + i} } \over {\left( {q_k  - 1} \right)^{i + 1} }}}
\Delta ^i 0^\alpha.
$$

Then, using theorem 2.1 and the identity (2.2), for $A_1$ we
obtain
    $$
A_1  = {{ - q_k^{N - s + 1}  + ( - 1)^\alpha  q_k^{s + 1} } \over
{\left( {q_k  - 1} \right)^{\alpha  + 1} }}E_{\alpha  - 1} (q_k
)\,.\eqno (3.6)
$$
For $A_2$ also, using theorem 2.1 and the identity (2.2) we have
$$
A_2  = {{ - q_k^{N - s + 1}  + ( - 1)^\alpha  q_k^{s + 1} } \over
{\left( {q_k  - 1} \right)^{\alpha  + 1} }}E_{\alpha  - 1} (q_k
)\,.\eqno (3.7)
$$
From equalities (3.6) and (3.7) we will obtain the statement of
lemma 3.2. Lemma 3.2 is proved.

\textbf{Proof of lemma 3.3.} Here we will use lemma 3.1 and the
equality (1.9). The optimal coefficients $C_\beta$ and the nodes
$x_\beta$ have to satisfy the equality (1.9), i.e.
$$
\sum\limits_{\beta=0}^NC_{\beta}x_{\beta}^{\alpha}=\frac{1}{\alpha+1},\ \ \
\alpha=0,1,2,...,m-1.
\eqno (3.8)
$$
From symmetry of the nodes (1.5) follows that
$$
C_\beta=C_{N - \beta }.  \eqno (3.9)
$$
Here the case $\alpha=0$  we consider separately.

Let $\alpha=0$ then from (3.8) we have
$$
\sum\limits_{\beta=0}^NC_{\beta}=1.\eqno (3.10)
$$
Taking into account (3.9) and using lemma 3.1 for the left side
of (3.10) we get
$$
\sum\limits_{\beta=0}^NC_{\beta}=2\sum\limits_{\beta=0}^{t-1}C_{\beta}+
\sum\limits_{\beta=t}^{N-t}h\left(1+\sum\limits_{k=1}^{m-1}d_k(q_k^{\beta}+
q_k^{N-\beta})\right)=
$$
$$
=2\sum\limits_{\beta=0}^{t-1}C_{\beta}-h\left(2t-1-\sum\limits_{k=1}^{m-1}
d_k\frac{q_k^t-q_k^{N-t+1}}{1-q_k}\right)+1.
$$

Hence, taking into account (3.10), we obtain
$$
\sum\limits_{\beta=0}^{t-1}C_{\beta}=h\left(\frac{2t-1}{2}-\sum\limits_{k=1}^{m-1}
d_k\frac{q_k^t-q_k^{N-t+1}}{1-q_k}\right).\eqno (3.11)
$$
This is the first equation of the system (3.3) corresponding to
the case $\alpha=0$. Thus we have proved the lemma for the case
$\alpha=0$.

Let, now, $\alpha=1,2,...,m-1$. Then, using symmetry of the nodes
(1.5) and keeping in mind the equality (3.9), for the left side of
(3.8) we obtain
$$
\sum\limits_{\beta  = 0}^N {C_\beta  x_\beta ^\alpha   = }
\sum\limits_{\beta  = 0}^{t - 1} {C_\beta  \left( {x_\beta ^\alpha
+ (1 - x_\beta ^{} )^\alpha  } \right) + \sum\limits_{\beta  =
t}^{N - t} {C_\beta  x_\beta ^\alpha   = Y_1  + Y_2 ,} }
$$
where
$$
Y_1  = \sum\limits_{\beta  = 0}^{t - 1} {C_\beta  \left( {x_\beta
^\alpha   + (1 - x_\beta ^{} )^\alpha  } \right),\,\,\,\,\,\,\,Y_2
= \sum\limits_{\beta  = t}^{N - t} {C_\beta  x_\beta ^\alpha  .} }
$$
Substituting the expression (3.1) of the optimal coefficients $C_\beta $
into $Y_2$, adding and subtracting the expressions
$$
\sum\limits_{\beta  = 1}^{t - 1} {h\left( {1 + \sum\limits_{k =
1}^{m - 1} {d_k \left( {q_k^\beta   + q_k^{N - \beta } } \right)}
} \right)(h\beta )^\alpha  \,}
$$
and
$$
\sum\limits_{\beta  = N - t + 1}^{N-1} {h\left( {1 + \sum\limits_{k =
1}^{m - 1} {d_k \left( {q_k^\beta   + q_k^{N - \beta } } \right)}
} \right)(h\beta )^\alpha  },
$$
we have
$$
Y_2  = h^{\alpha  + 1} \left( {\sum\limits_{\beta  = 1}^{N - 1}
{\beta ^\alpha   + \sum\limits_{k = 1}^{m - 1} {d_k \left(
{\sum\limits_{\beta  = 1}^{N - 1} {q_k^\beta  \beta ^\alpha   +
\sum\limits_{\beta  = 1}^{N - 1} {q_k^{N - \beta } \beta ^\alpha
} } } \right)} } } \right) -
$$
$$
-h\sum\limits_{\beta  = 1}^{t - 1} {\left( {\left( {h\beta }
\right)^\alpha   + \left( {1 - h\beta } \right)^\alpha  }
\right)\left( {1 + \sum\limits_{k = 1}^{m - 1} {d_k \left(
{q_k^\beta   + q_k^{N - \beta } } \right)} } \right).}
$$
Using equalities (2.4), (2.5), (2.6) and binomial formula, then
grouping in powers of $h$ and taking into account lemma 3.2, for
$Y_2$ we have
$$
Y_2  = h^{\alpha  + 1} \left[ {\sum\limits_{k = 1}^{m - 1} {d_k
\sum\limits_{i = 1}^\alpha  {{{ - q_k^{N + i}  + ( - 1)^i q_k }
\over {(1 - q_k^{} )^{i + 1} }}\Delta ^i 0^\alpha   -
\sum\limits_{\beta  = 1}^{t - 1} {\beta ^\alpha  \left( {1 +
\sum\limits_{k = 1}^{m - 1} {d_k \left( {q_k^\beta   + q_k^{N -
\beta } } \right)} } \right)} } } } \right](1 +
(-1)^\alpha)+
$$
$$
+ \sum\limits_{j =1}^{\alpha  - 1} {{{\alpha!h^{j + 1} } \over
{j!(\alpha  - j)!}}\left[ {{{B_{j + 1} } \over {j + 1}} +
\sum\limits_{k = 1}^{m - 1} {d_k \sum\limits_{i = 1}^j {{{ -
q_k^{N + i}+(-1)^i q_k }\over{(1 - q_k )^{i + 1} }}\Delta
^i 0^j +}   } } \right.}
$$
$$
+\left.( - 1)^{j + 1} \sum\limits_{\beta  = 1}^{t - 1} {\beta ^j \left(
{1 + \sum\limits_{k = 1}^{m - 1} {d_k \left( {q_k^\beta   + q_k^{N
- \beta } } \right)} } \right)}\right]-
$$
$$
-h\left(\frac{2t-1}{2}-\sum\limits_{k=1}^{m-1}d_k
\frac{q_k^t-q_k^{N-t+1}}{1-q_k}\right)
 + {1 \over {\alpha  + 1}}.
$$
Now we consider $Y_1$. Taking into account that when
$\beta=\overline{0,t-1}$ then from (1.5) the nodes
$x_{\beta}=\eta_{\beta}h$ and applying binomial formula, after
that, grouping in powers of $h$, for $Y_1$  we have
$$
Y_1  = \sum\limits_{\beta  = 0}^{t - 1} {C_\beta  (x_\beta ^\alpha
+ (1 - x_\beta  )^\alpha  ) = } \sum\limits_{\beta  = 0}^{t - 1}
{C_\beta  \left( {\left( {\eta _\beta  h} \right)^\alpha   +
\sum\limits_{j = 0}^\alpha  {{{\alpha !( - \eta _\beta  h)^j }
\over {j!(\alpha  - j)!}}} } \right) = }
$$
$$
= h^{\alpha  + 1} \sum\limits_{\beta  = 0}^{t - 1} {h^{ - 1}
C_\beta  \eta _\beta ^\alpha  (1 + ( - 1)^\alpha  ) +
\sum\limits_{j = 0}^{\alpha  - 1} {{{\alpha !h^{j + 1} } \over
{j!(\alpha  - j)!}}\sum\limits_{\beta  = 0}^{t - 1} {h^{ - 1}
C_\beta  ( - \eta _\beta  )^j .} } }
$$
Hence using (3.11) we obtain
$$
Y_1=
h^{\alpha  + 1} \sum\limits_{\beta  = 0}^{t - 1} h^{ - 1}
C_\beta  \eta _\beta ^\alpha  (1 + ( - 1)^\alpha  ) +
$$
$$
+\sum\limits_{j = 1}^{\alpha  - 1} {{{\alpha !h^{j + 1} } \over
{j!(\alpha  - j)!}}\sum\limits_{\beta  = 0}^{t - 1} {h^{ - 1}
C_\beta  ( - \eta _\beta  )^j } } +h\left(\frac{2t-1}{2}-\sum\limits_{k=1}^{m-1}d_k
\frac{q_k^t-q_k^{N-t+1}}{1-q_k}\right).
$$
Now, adding $Y_1$, $Y_2$ and the result substituting to the left
side of the equality (3.8), we obtain
$$
h^{\alpha  + 1} \left[
\sum\limits_{k =1}^{m - 1}d_k \sum\limits_{i = 1}^\alpha
\frac{- q_k^{N + i}+(-1)^i q_k}{(1 - q_k )^{i + 1}}\Delta ^i 0^\alpha- \right.
$$
$$
-\sum\limits_{\beta= 1}^{t - 1}\beta^\alpha
\left(1+\sum\limits_{k = 1}^{m - 1} {d_k \left( {q_k^\beta   + q_k^{N -
\beta } } \right)} \right) +
\left.\sum\limits_{\beta  = 0}^{t - 1}h^{ - 1}
C_\beta  \eta _\beta ^\alpha   \right]\left({1 +( - 1)^\alpha
}\right)+
$$
$$
+ \sum\limits_{j = 1}^{\alpha  - 1} {{\alpha!h^{j + 1} } \over
{j!(\alpha  - j)!}}\left[{{B_{j + 1}}\over {j + 1}} +
\sum\limits_{k = 1}^{m - 1}d_k \sum\limits_{i = 1}^j
\frac{-q_k^{N + i}+(-1)^i q_k }{(1 - q_k )^{i + 1}}\Delta
^i 0^j  +  \right.
$$
$$
+(-1)^{j+1}\left.\sum\limits_{\beta  = 1}^{t - 1} \beta
^j \left(1 + \sum\limits_{k = 1}^{m - 1}d_k\left(q_k^\beta+
q_k^{N - \beta } \right)\right) + \sum\limits_{\beta  =
0}^{t - 1} {h^{ - 1} C_\beta(-\eta_\beta)^j }\right]
= 0. \eqno (3.12)
$$
Clearly, that the left side of (3.12) is the polynomial of degree
$\alpha+1$ with respect to $h$. Since this polynomial is equal to
zero then all coefficients of the polynomial are zero.

When $\alpha$ is odd from (3.12) we get system of equations
$$
\sum\limits_{k=1}^{m-1}d_k\left\{
\sum\limits_{i=1}^j\frac{-q_k^{N+i}+(-1)^iq_k}{(1-q_k)^{i+1}}\Delta^i0^j+
(-1)^{j+1}\sum\limits_{\beta=1}^{t-1}\beta^j
\left(q_k^{\beta}+q_k^{N-\beta}\right)\right\}+
$$
$$
+\frac{B_{j+1}}{j+1}+(-1)^j\sum\limits_{\beta=0}^{t-1}\left(C_{\beta}h^{-1}
\eta_\beta^j-\beta^j\right)=0, \ \ \ \ j=1,2,...,\alpha-1,
$$
which is the part of the system (4.22), because $1\le \alpha\le
m-1$.

When $\alpha$ is even from (3.12) we get following system of
equations
$$
\sum\limits_{k=1}^{m-1}d_k\left\{
\sum\limits_{i=1}^j\frac{-q_k^{N+i}+(-1)^iq_k}{(1-q_k)^{i+1}}\Delta^i0^j+
(-1)^{j+1}\sum\limits_{\beta=1}^{t-1}\beta^j
\left(q_k^{\beta}+q_k^{N-\beta}\right)\right\}+
$$
$$
+\frac{B_{j+1}}{j+1}+(-1)^j\sum\limits_{\beta=0}^{t-1}\left(C_{\beta}h^{-1}
\eta_\beta^j-\beta^j\right)=0, \ \ \ \ j=1,2,...,\alpha-1,
$$
which is the part of the system (4.22) and one new equation for
each even $\alpha$
$$
\sum\limits_{k=1}^{m-1}d_k\left[\sum\limits_{i=1}^\alpha\frac{-q_k^{N+i}+(-1)^iq_k}{(1-q_k)^{i+1}}\Delta^i0^\alpha-
\sum\limits_{\beta=1}^{t-1}\beta^\alpha\left(q_k^{\beta}+q_k^{N-\beta}\right)\right]+
\sum\limits_{\beta=0}^{t-1}(C_\beta
h^{-1}\eta_\beta^\alpha-\beta^{\alpha})=0.
$$
Hence and from (3.11) for the optimal coefficients $C_{\beta}$,
$\beta=0,...,t-1$ we will get following system of equations
$$
\sum\limits_{\beta=0}^{t-1}C_{\beta}=h\left(\frac{2t-1}{2}-\sum\limits_{k=1}^{m-1}
d_k\frac{q_k^t-q_k^{N-t+1}}{1-q_k}\right),\eqno (3.13)
$$
$$
\sum\limits_{\beta=0}^{t - 1}C_\beta  \eta _\beta ^\alpha   =
h\left[ \sum\limits_{\beta  = 1}^{t - 1} \beta ^\alpha -
\sum\limits_{k = 1}^{m - 1} d_k \left( \sum\limits_{i = 1}^\alpha
\frac{ - q_k^{N + i} + ( - 1)^i q_k }{(1 - q_k )^{i + 1} }\Delta
^i 0^\alpha-\sum\limits_{\beta  = 1}^{t - 1} {\beta ^\alpha\left(
{q_k^\beta+ q_k^{N - \beta } } \right)} \right)\right],\eqno
(3.14)
$$
where $\alpha=\left\{
\begin{array}{ll}
2,4,...,m-2 &\mbox{ when }m-\mbox{ even},\\
2,4,...,m-1 &\mbox{ when }m-\mbox{ odd}.\\
\end{array}
\right. $\\
Thus for solvability of the system (3.13), (3.14) (with respect to
$C_{\beta}$, $\beta=\overline{0,t-1}$) $t$ have to take on a value
$\frac{m}{2}$ when $m$ is even and $\left[\frac{m}{2}\right]+1$
when $m$ is odd, where $[a]$ is the integer part of the number
$a$.

Now, let $Y_3  = \sum\limits_{\beta  = 1}^{t - 1} {\beta ^\alpha
\left( {q_k^\beta   + q_k^{N - \beta } } \right)}$. Then, using
formula (2.4) and lemma 3.2 taking into account that $\alpha$ is
even natural number,  we have
$$
Y_3  = \sum\limits_{i = 1}^\alpha  {{{ - q_k^{N + i}  + ( - 1)^i
q_k } \over {(1 - q_k )^{i + 1} }}\Delta ^i 0^\alpha   +
\sum\limits_{i = 0}^\alpha  {{{( - 1)^i q_k^{t + i}  - q_k^{N - t
+ 1} } \over {(q_k  - 1)^{i + 1} }}\Delta ^i t^\alpha } }.\eqno
(3.15)
$$
Putting the equality (3.15) into (3.14), taking into account
(3.13) we will get the
assertion of lemma 3.3. Lemma 3.3 is proved.\\[0.2cm]

\textbf{Proof of lemma 4.} Using formula of the binomial
coefficients we will get
$$
C_m^l C_l^k  = {{m!} \over {l! \cdot (m - l)!}}{{l!} \over {k!
\cdot (l - k)!}} =
$$
$$
= {{m!} \over {k! \cdot (m - k)!}} \cdot {{(m - k)!} \over {(l -
k)! \cdot \left( {\left( {m - k} \right) - \left( {l - k} \right)}
\right)!}} = C_m^k C_{m - k}^{l - k}.\eqno (3.16)
$$

Applying the equality (3.16) to the left side of (3.4) we obtain
    $$
C_m^{m - 1} C_{m - 1}^n  - C_m^{m - 2} C_{m - 2}^n  + ... + ( -
1)^{m - n - 2} C_m^{n + 1} C_{n + 1}^n  + ( - 1)^{m - n - 1} C_m^n
C_n^n  =
$$
$$
= C_m^n C_{m - n}^{m - n - 1}  - C_m^n C_{m - n}^{m - n - 2}  +
... + ( - 1)^{m - n - 1} C_m^n C_{m - n}^0  =
$$
$$
= C_m^n \left( {C_{m - n}^{m - n - 1}  - C_{m - n}^{m - n - 2}  +
... + ( - 1)^{m - n - 1} C_{m - n}^0 } \right) =
$$
$$
= C_m^n \left( {C_{m - n}^{m - n}  - \left( {C_{m - n}^{m - n}  -
C_{m - n}^{m - n - 1}  + ... + ( - 1)^{m - n} C_{m - n}^0 }
\right)} \right) =
$$
$$
=
 C_m^n \left( {C_{m - n}^{m - n}  - (1 - 1)^{m -
n} } \right) = C_m^n,
$$
Lemma 3.4 is proved.

\textbf{Proof of lemma 3.5.} The assertion of this lemma we will
prove by induction.\\ For $k = 1$:
$$
y = \sum\limits_{i = 1}^1 {{{y^{\left[ i \right]} } \over
{i!}}\Delta ^i 0^1  = } {{y^{\left[ 1 \right]} } \over {1!}}\Delta
^1 0^1  = y,
$$
For $k = 2$:
$$
y^2  = \sum\limits_{i = 1}^2 {{{y^{\left[ i \right]} } \over
{i!}}\Delta ^i 0^2  = } {{y^{\left[ 1 \right]} } \over {1!}}\Delta
^1 0^2  + {{y^{\left[ 2 \right]} } \over {2!}}\Delta ^2 0^2  = y +
y(y - 1) = y^2 .
$$
Now, suppose that the equality (3.5) is true for $k = \alpha$,
i.e.
    $$
y^\alpha   = \sum\limits_{i = 1}^\alpha  {{{y^{\left[ i \right]} }
\over {i!}}\Delta ^i 0^\alpha   = } {{y^{\left[ 1 \right]} } \over
{1!}}\Delta ^1 0^\alpha   + {{y^{\left[ 2 \right]} } \over
{2!}}\Delta ^2 0^\alpha   + ... + {{y^{\left[ \alpha  \right]} }
\over {\alpha !}}\Delta ^\alpha  0^\alpha.\eqno (3.17)
$$
We will show that the equality (3.5) is also true for the case
$k=\alpha  + 1$. For this we will use following equalities
    $$
y^{\left[ {\alpha  + 1} \right]}  = y^{\left[ \alpha  \right]} (y
- \alpha ),\eqno (3.18)
$$
    $$
\Delta ^i 0^{\alpha  + 1}  = i\left( {\Delta ^{i - 1} 0^\alpha   +
\Delta ^i 0^\alpha  } \right).\eqno (3.19)
$$
Multiplying both parts of the equality (3.17) by $y$ and keeping
in mind (3.18), (3.19), we obtain
    $$
y^{\alpha  + 1}  = \sum\limits_{i = 1}^{\alpha  + 1} {{{y^{\left[
i \right]} } \over {i!}}\Delta ^i 0^{\alpha  + 1}  = } {{y^{\left[
1 \right]} } \over {1!}}\Delta ^1 0^{\alpha  + 1}  + {{y^{\left[ 2
\right]} } \over {2!}}\Delta ^2 0^{\alpha  + 1}  + ... +
{{y^{\left[ {\alpha  + 1} \right]} } \over {(\alpha  + 1)!}}\Delta
^{\alpha  + 1} 0^{\alpha  + 1} .
$$
Thus, for any $k \in N$ the equality (3.5) is true. Lemma 3.5 is
proved.

\section{The main results}

\subsection{Coefficients of optimal quadrature formulas}

For the coefficients of optimal quadrature formulas of the form
(1.1) following theorem holds.

\textbf{Theorem 4.1.} {\it The coefficients of optimal quadrature
formulas of the form (1.1) with the error functional (1.2) and the
nodes (1.5) in Sobolev space $L_2^{(m)}(0,1)$ are expressed by
formula
$$
C_\beta   = h\left( {1 + \sum\limits_{k = 1}^{m - 1} {d_k \left(
{q_k^\beta   + q_k^{N - \beta } } \right)} } \right),\,\,\,\,\,\,t
\le \beta  \le N - t,\eqno (4.1)
$$
where $d_k$, $(k=\overline{1,m-1})$ satisfy following system of
$m-1$ linear equations:
$$
\sum\limits_{k=1}^{m-1}d_k\sum\limits_{i=1}^j
\frac{-q_k^{t+1}+(-1)^i q_k^{N-t+i}}{(q_k-1)^{i+1}}\Delta^i0^j=
\frac{t^{j+1}-B_{j+1}}{j+1}-\sum\limits_{\beta=0}^{t-1}C_{\beta}h^{-1}
(t-\eta_{\beta})^j, \eqno (4.2)
$$
$$
j=1,2,3,...,m-1,
$$
here the coefficients $C_\beta   = C_{N - \beta }$ $(\beta
=0,1,...,t - 1)$ are determined from the system (3.3), $q_k$ are
roots of Euler-Frobenius polynomial $E_{2m -2}(q)$, $|q_k |< 1$.}

\textbf{Proof.} As said before the optimal coefficients
$C_{\beta}$ ($\beta=\overline{0,N}$) are solution of the system
(1.8)-(1.9). In the lemma 3.1 we have obtained representation of
the optimal coefficients $C_{\beta}$ for $\beta=\overline{t,N-t}$.
In the proof of the lemma 3.3, using the result of lemma 3.1, we
have obtained the system (3.3) for the optimal coefficients
$C_{\beta}$ ($\beta=\overline{0,t-1}$). So we conclude that in
order to solve the system (1.8)-(1.9) it is sufficient to find
unknown parameters $d_k$ and unknown polynomial $P_{m - 1}
(x_\beta) = \sum\limits_{k = 1}^{m - 1} {\lambda _\alpha x_\beta
^\alpha}$. For this in (1.8) instead of $C_{\beta}$,
$\beta=\overline{t,N-t}$ substituting the expression (3.1) we will
get polynomials of degree $2m$ with respect to $x_{\beta}$ in both
sides of (1.8). Equating coefficients of same powers of $x_\beta$
we will get the system of $m-1$ linear equations for unknowns
$d_k$ and we will find the coefficients of unknown polynomial
$P_{m-1}(x_{\beta})$. Thus the proof will be complete.

Further we will give detailed explanation of the proof of the
theorem.

Now we consider the equality (1.8)
$$
\sum\limits_{\gamma = 0}^N {C_\gamma{{|x_\beta   - x_\gamma |^{2m
- 1} } \over {2 \cdot (2m - 1)!}} + P_{m - 1} (x_\beta  ) =
{{x_\beta ^{2m}  + (1 - x_\beta  )^{2m}} \over {2\left({2m}
\right)!}}\ \ \ \ \mbox{as}\ \ \ x_\beta \in [0,1].}\eqno (4.3)
$$
We denote by
$$
g(x_\beta) = \sum\limits_{\gamma  = 0}^N {C_\gamma  {{|x_\beta -
x_\gamma  |^{2m - 1} } \over {2 \cdot (2m - 1)!}},} \eqno (4.4)
$$
 $$
f(x_\beta  ) = {{x_\beta ^{2m}  + (1 - x_\beta  )^{2m} } \over
{2\left( {2m} \right)!}}.\eqno (4.5)
$$
Now we consider the function $g(x_\beta)$ and let $t \le \beta \le
N - t$. Then
$$ g(x_\beta  ) = \sum\limits_{\gamma  = 0}^{t -
1} {C_\gamma {{(x_\beta   - x_\gamma  )^{2m - 1} } \over {(2m -
1)!}} + } \sum\limits_{\gamma  = t}^\beta  {C_\gamma  {{(x_\beta -
x_\gamma  )^{2m - 1} } \over {(2m - 1)!}} - } \sum\limits_{\gamma
= 0}^N {C_\gamma  {{(x_\beta   - x_\gamma )^{2m - 1} } \over {2(2m
- 1)!}}.}\eqno (4.6)
$$
Further we denote by
$$
\psi _1 (x_\beta  ) = \sum\limits_{\gamma  = 0}^{t - 1} {C_\gamma
{{(x_\beta   - x_\gamma  )^{2m - 1} } \over {(2m - 1)!}},}\eqno
(4.7)
$$
$$
\psi _2 (x_\beta  ) = \sum\limits_{\gamma  = t}^\beta  {C_\gamma
{{(x_\beta   - x_\gamma  )^{2m - 1} } \over {(2m - 1)!}}}, \eqno
(4.8)
$$
$$
\psi _3 (x_\beta  ) =  - \sum\limits_{\gamma  = 0}^N {C_\gamma
{{(x_\beta   - x_\gamma  )^{2m - 1} } \over {2(2m - 1)!}}.}\eqno
(4.9)
$$
In the equality (4.8), using (3.1), adding and subtracting
following expression
$$
{{h^{2m} } \over {(2m - 1)!}}\left( {\sum\limits_{\gamma  = 1}^{t
- 1} {\left( {\beta  - \gamma } \right)^{2m - 1}  + \sum\limits_{k
= 1}^{m - 1} {d_k \left( {\sum\limits_{\gamma  = 1}^{t - 1}
{q_k^\gamma  \left( {\beta  - \gamma } \right)^{2m - 1}  +
\sum\limits_{\gamma  = 1}^{t - 1} {q_k^{N - \gamma } \left( {\beta
- \gamma } \right)^{2m - 1} } } } \right)} } } \right),
$$
we get
$$
\psi _2 (x_\beta  ) = {{h^{2m} } \over {(2m - 1)!}}\left[
{\sum\limits_{\gamma  = 1}^\beta  {\left( {\beta  - \gamma }
\right)^{2m - 1}  + \sum\limits_{k = 1}^{m - 1} {d_k \left(
{\sum\limits_{\gamma  = 1}^\beta  {q_k^\gamma  \left( {\beta  -
\gamma } \right)^{2m - 1}  + \sum\limits_{\gamma  = 1}^\beta
{q_k^{N - \gamma } \left( {\beta  - \gamma } \right)^{2m - 1} } }
} \right)} } } \right] -
$$
$$
- {{h^{2m} } \over {(2m - 1)!}}\left[ {\sum\limits_{\gamma  =
1}^{t - 1} {\left( {\beta  - \gamma } \right)^{2m - 1}  +
\sum\limits_{k = 1}^{m - 1} {d_k \left( {\sum\limits_{\gamma  =
1}^{t - 1} {q_k^\gamma  \left( {\beta  - \gamma } \right)^{2m - 1}
+ \sum\limits_{\gamma  = 1}^{t - 1} {q_k^{N - \gamma } \left(
{\beta  - \gamma } \right)^{2m - 1} } } } \right)} } } \right].
$$
Hence replacing  $\beta - \gamma$ with $\gamma$ and using
equalities (2.4), (2.5)
$$
\psi_2(x_{\beta})=\frac{h^{2m}}{(2m-1)!}\Bigg[
\sum\limits_{j=1}^{2m}\frac{(2m-1)!B_{2m-j}}{j!\ (2m-j)!}\beta^j+
$$
$$
+\sum\limits_{k=1}^{m-1}\Bigg\{q_k^{\beta}\left( \frac{q_k}{q_k-1}
\sum\limits_{i=0}^{2m-1}\frac{\Delta^i0^{2m-1}}{(q_k-1)^i}-
\frac{q_k^{1-\beta}}{q_k-1}\sum\limits_{i=0}^{2m-1}
\frac{\Delta^i\beta^{2m-1}}{(q_k-1)^i}\right)+
$$
$$
+q_k^{N-\beta}\left( \frac{1}{1-q_k}
\sum\limits_{i=0}^{2m-1}\left(\frac{q_k}{1-q_k}\right)^i\Delta^i0^{2m-1}-
\frac{q_k^{\beta}}{1-q_k}\sum\limits_{i=0}^{2m-1}
\left(\frac{q_k}{1-q_k}\right)^i\Delta^i\beta^{2m-1}\right)\Bigg\}+
$$
$$
-\sum\limits_{\gamma=1}^{t-1}(\beta-\gamma)^{2m-1}\left(1+\sum\limits_{k=1}
^{m-1}d_k\left(q_k^{\gamma}+q_k^{N-\gamma}\right)\right)\Bigg].
$$
Taking into account that $q_k$ are roots of the Euler-Frobenius
polynomial $E_{2m-2}(q)$ and  using theorem 2.1, for
$\psi_2(x_\beta)$ we have
$$
\psi_2(x_{\beta})=\frac{h^{2m}}{(2m-1)!}\Bigg[\sum\limits_{j=1}^{2m}
\frac{(2m-1)!B_{2m-j}}{j!\ (2m-j)!}\beta^j+\sum\limits_{k=1}^{m-1}
d_k\sum\limits_{i=0}^{2m-1}\frac{-q_k+(-1)^iq_k^{N+i}}{(q_k-1)^{i+1}}
\Delta^i\beta^{2m-1}-
$$
$$
-\sum\limits_{\gamma=1}^{t-1}(\beta-\gamma)^{2m-1}\left(1+\sum\limits_{k=1}
^{m-1}d_k\left(q_k^{\gamma}+q_k^{N-\gamma}\right)\right)\Bigg].
\eqno (4.10)
$$
Finally, using binomial formula and (2.6) from (4.10) we get
$$
\psi _2 (x_\beta  ) = {{(h\beta) ^{2m} } \over {(2m)!}} + {{h^{2m}
\beta ^{2m - 1} } \over {(2m - 1)!}}\left[ -\frac{2t-1}{2}+
\sum\limits_{k=1}^{m-1}d_k\frac{q_k^t-q_k^{N-t+1}}{1-q_k}\right] +
$$
$$
+{{h^{2m} } \over {(2m - 1)!}}\sum\limits_{j = 1}^{2m - 2} {{{(2m
- 1)!B_{j + 1} \beta ^{2m - 1 - j} } \over {(j + 1)!(2m - 1 -
j)!}} + } {{h^{2m} } \over {(2m - 1)!}}\sum\limits_{j = 1}^{2m -
1} {{{(2m - 1)!\beta ^{2m - 1 - j} } \over {j!(2m - 1 - j)!}}
\times }
$$
$$
\times \left[ {\sum\limits_{k = 1}^{m - 1} {d_k \sum\limits_{i =
1}^j {{{ - q_k  + ( - 1)^i q_k^{N + i} } \over {(q_k  - 1)^{i + 1}
}}\Delta ^i 0^j  - \sum\limits_{\gamma  = 1}^{t - 1} {( - \gamma
)^j \left( {1 + \sum\limits_{k = 1}^{m - 1} {d_k \left(
{q_k^\gamma   + q_k^{N - \gamma } } \right)} } \right)} } } }
\right].\eqno (4.11)
$$
Now we consider the equality (4.7). Using binomial formula,
keeping in mind (1.5) and (3.3) when $\alpha  = 0$, for
$\psi_1(x_\beta )$ we have
$$
\psi _1 (x_\beta  )
=\sum\limits_{\gamma=0}^{t-1}C_{\gamma}\frac{(x_{\beta}-x_{\gamma})^{2m-1}}{(2m-1)!}=
\sum\limits_{\gamma=0}^{t-1}C_{\gamma}\frac{(h\beta-\eta_\gamma
h)^{2m-1}}{(2m-1)!}=
$$
$$
={{h^{2m} } \over {(2m - 1)!}}\sum\limits_{j = 1}^{2m - 1} {{{(2m
- 1)!\beta ^{2m - 1 - j} } \over {j!(2m - 1 -
j)!}}\sum\limits_{\gamma  = 0}^{t - 1} {C_\gamma  h^{ - 1} ( -
\eta _\gamma  )^j  + } }
$$
$$
+ {{h^{2m} \beta ^{2m - 1} } \over {(2m - 1)!}}\left[ {{{2t - 1}
\over 2} - \sum\limits_{k = 1}^{m - 1} {d_k {{q_k^t  - q_k^{N - t
+ 1} } \over {1 - q_k }}} } \right].\eqno (4.12)
$$

Further using binomial formula and the equality (1.9) from (4.9)
we obtain
$$
\psi _3 (x_\beta  ) =  - \sum\limits_{j = 0}^{m - 1} {{{( - 1)^j
x_\beta ^{2m - 1 - j} } \over {2(j + 1)!(2m - 1 - j)!}} -
\sum\limits_{j = 0}^{m - 1} {{{( - 1)^{m + j} x_\beta ^{m - 1 - j}
} \over {2(m + j)!(m - 1 - j)!}}\sum\limits_{\gamma  = 0}^N
{C_\gamma  x_\gamma ^{m + j} .} } }\eqno (4.13)
$$
Similarly, using binomial formula in the equality (4.5) for
$f_m(x_{\beta})$ we have
$$
f_m (x_\beta  ) = {{x_\beta ^{2m} } \over {(2m)!}} +
\sum\limits_{j = 0}^{m - 1} {{{( - 1)^{j + 1} x_\beta ^{2m - 1 -
j} } \over {2(j + 1)!(2m - 1 - j)!}} + \sum\limits_{j = 0}^{m - 1}
{{{( - 1)^{m + j} x_\beta ^{m + 1 + j} } \over {2(m + 1 + j)!(m -
1 - j)!}}.} }\eqno (4.14)
$$
Substituting equalities (4.11), (4.12), (4.13), (4.14)  to  (4.3)
and after calculations we obtain
$$
{{h^{2m} } \over {(2m - 1)!}}\sum\limits_{j = 1}^{2m - 2} {{{(2m -
1)!B_{j + 1} \beta ^{2m - 1 - j} } \over {(j + 1)!(2m - 1 - j)!}}
+ {{h^{2m} } \over {(2m - 1)!}}\sum\limits_{j = 1}^{2m - 1} {{{(2m
- 1)!\beta ^{2m - 1 - j} } \over {j!(2m - 1 - j)!}} \times } }
$$
$$
\times \left[ {( - 1)^j \sum\limits_{\gamma  = 0}^{t - 1} {\left(
{C_\gamma  h^{ - 1} \eta _\gamma ^j  - \gamma ^j } \right) +
\sum\limits_{k = 1}^{m - 1} {d_k \left( {\sum\limits_{i = 1}^j {{{
- q_k  + ( - 1)^i q_k^{N + i} } \over {(q_k  - 1)^{i + 1} }}\Delta
^i 0^j  - \sum\limits_{\gamma  = 1}^{t - 1} {( - \gamma )^j \left(
{q_k^\gamma   + q_k^{N - \gamma } } \right)} } } \right)} } }
\right] +
$$
$$
+ \sum\limits_{j = 0}^{m - 1} {{{( - 1)^{m - j} x_\beta ^{m - 1 -
j} } \over {2(m + j)!(m - 1 - j)!}}\left[ {{1 \over {m + 1 + j}} -
\sum\limits_{\gamma  = 0}^N {C_\gamma  x_\gamma ^{m + j} } }
\right]}  =  - P_{m - 1} (x_\beta  ).\eqno (4.15)
$$
Keeping in mind the equality (4.13) and designations (4.4), (4.5),
from (4.15) one can see that the difference $g(x_\beta  ) -
f_m(x_\beta )$ is polynomial of degree $2m - 2$ with respect to
$x_\beta$, i.e.
$$
g(x_\beta) - f_m(x_\beta) = \sum\limits_{j = 0}^{2m - 2} {a_j
x_\beta ^j ,\,\,\,\,\,\,t \le \beta  \le N - t.}\eqno (4.16)
$$
Here
$$
a_j  = \left\{
\begin{array}{l}
{b_j \,\,\,\,\,\mbox{as}\,\,\,\,m \le j\le 2m - 2,} \\
{b_j  + {{( - 1)^{j + 1} } \over {2j!(2m - j - 1)!}}\left[ {{1
\over {2m - j}} - \sum\limits_{\gamma  = 0}^N {C_\gamma  x_\gamma
^{2m - j - 1} } } \right]\,\,\,\,\mbox{as}\,\,\,\,1 \le j \le m - 1,} \\
{{{h^{2m} } \over {(2m - 1)!}}\left[ {( - 1)^{2m - 1}
\sum\limits_{\gamma  = 0}^{t - 1} {\left( {C_\gamma  h^{ - 1} \eta
_\gamma ^{2m - 1}  - \gamma ^{2m - 1} } \right) + } } \right.} \\
{ + \sum\limits_{k = 1}^{m - 1}d_k {\left( {\sum\limits_{i =
1}^{2m - 1} {{{ - q_k  + ( - 1)^i q_k^{N + i} } \over {(q_k  -
1)^{i + 1} }}\Delta ^i 0^{2m - 1}  + \sum\limits_{\gamma  = 1}^{t
- 1} {\gamma ^{2m - 1} \left( {q_k^\gamma   + q_k^{N - \gamma } }
\right)} } } \right) - } } \\
{ - {1 \over {2(2m - 1)!}}\left( {{1 \over {2m}} -
\sum\limits_{\gamma  = 0}^N {C_\gamma  x_\gamma ^{2m - 1} } }
\right)\,\,\,\mbox{as}\,\,j = 0,}
 \end{array}
\right.\eqno (4.17)
$$
where
$$
b_j  = {{h^{2m - j} } \over {j!(2m - j - 1)!}}\left[ {{{B_{2m - j}
} \over {2m - j}} + ( - 1)^{2m - j - 1} \sum\limits_{\gamma  =
0}^{t - 1} {\left( {C_\gamma  h^{ - 1} \eta _\gamma ^{2m - j - 1}
- \gamma _{}^{2m - j - 1} } \right) + } } \right.
$$
$$
+ \sum\limits_{k = 1}^{m - 1} {d_k \left( {\sum\limits_{i = 1}^{2m
- j - 1} {{{ - q_k  + ( - 1)^i q_k^{N + i} } \over {(q_k  - 1)^{i
+ 1} }}\Delta ^i 0^{2m - 1 - j}  + ( - 1)^{2m - j}
\sum\limits_{\gamma  = 1}^{t - 1} {\gamma ^{2m - j - 1} \left(
{q_k^\gamma   + q_k^{N - \gamma } } \right)} } } \right).}
$$
On the other hand from (4.15) we obtain
$$
g(x_\beta  ) - f_m (x_\beta  )=-P_{m - 1}(x_\beta).\eqno (4.18)
$$
This equality takes placed if $b_j  = 0 $, as $m \le j \le 2m - 2$
or
$$
\sum\limits_{k = 1}^{m - 1} {d_k \left\{ {\sum\limits_{i = 1}^j
{{{ - q_k  + ( - 1)^i q_k^{N + i} } \over {(q_k  - 1)^{i + 1}
}}\Delta ^i 0^j  - ( - 1)^j \sum\limits_{\gamma  = 1}^{t - 1}
{\gamma ^j \left( {q_k^\gamma   + q_k^{N - \gamma } } \right)} } }
\right\} = }
$$
$$
= - {{B_{j + 1} } \over {j + 1}} - ( - 1)^j \sum\limits_{\gamma =
0}^{t - 1} {\left( {C_\gamma  h^{ - 1} \eta _\gamma ^j  - \gamma
^j } \right),\,\,\,\,j = 1,2,...,m - 1.}\eqno (4.19)
$$
From equalities (4.16) and (4.18) we will find unknown polynomial
$P_{m - 1}(x_\beta)$ of the system (1.8) - (1.9)
$$
P_{m - 1} (x_\beta  ) =  - \sum\limits_{j = 0}^{m - 1} {a_j
x_\beta ^j .}\eqno (4.20)
$$
Later, applying to the sum $A = \sum\limits_{\gamma  = 1}^{t - 1}
{\gamma ^j \left( {q_k^\gamma   + q_k^{N - \gamma } } \right)}$
formulas (2.4), (2.5), (2.6) and (3.2), we obtain
$$
A = ( - 1)^j \sum\limits_{i = 1}^j {{{ - q_k  + ( - 1)^i q_k^{N +
i} } \over {(q_k  - 1)^{i + 1} }}\Delta ^i 0^j  + } \sum\limits_{i
= 0}^j {{{ - q_k^{N - t + 1}  + ( - 1)^i q_k^{t + i} } \over {(q_k
- 1)^{i + 1} }}\sum\limits_{p = 0}^j {C_j^p } \Delta ^i 0^p t^{j -
p} .}\eqno (4.21)
$$
Substituting the equality (4.21) to the (4.19), after some
simplifications we have
$$
\sum\limits_{k = 1}^{m - 1} {d_k \left\{ {( - 1)^{j + 1}
\sum\limits_{p = 0}^j {C_j^p t^{j - p} } \sum\limits_{i = 0}^p {{{
- q_k^{N - t + 1}  + ( - 1)^i q_k^{t + i} } \over {(q_k  - 1)^{i +
1} }}\Delta ^i 0^p } } \right\} = }
$$
$$
=-{{B_{j + 1} } \over {j + 1}} - ( - 1)^j \sum\limits_{\gamma =
0}^{t - 1} {\left( {C_\gamma  h^{ - 1} \eta _\gamma ^j  - \gamma
^j } \right),\,\,\,\,j = 1,2,...,m - 1.}\eqno (4.22)
$$
Multiplying the first equation of the system (4.22) by $( - 1)^2
C_2^1 t$, adding with the second one, then, multiplying the first
equation by $( - 1)^3 C_3^1 t^2$, the second by $( - 1)^2 C_3^2 t$
of the system (4.22), adding with the third and  so on, continuing
by this way, and also taking into account lemmas 3.2, 3.3, 3.4,
3.5 and binomial formula for unknown parameters $d_k$ we get
following linear system of equations
$$
\sum\limits_{k=1}^Nd_k\sum\limits_{i=1}^j\frac{-q_k^{t+1}+(-1)^i
q_k^{N-t+i}}{(q_k-1)^{i+1}}\Delta^i0^j=
\frac{t^{j+1}-B_{j+1}}{j+1}- \sum\limits_{\beta=0}^{t-1}C_\beta
h^{-1}(t-\eta_\beta)^j,
$$
$$
j=1,2,...,m-1.
$$
The last system is the system (4.2) for $d_k$. This completes the
proof of the theorem 4.1.

\subsection{Norm of the error functional of the optimal quadrature formula}

In this section square of norm of the error functional (1.2) of
optimal quadrature formulas of the form (1.1) with the nodes (1.5)
is calculated.

Following is valid

\textbf{Theorem 4.2.} {\it Square of norm of the error functional
(1.2) of optimal quadrature formulas of the form (1.2) with the
nodes (1.5) on the space $L_2^{(m)} (0,1)$ have following form
$$
\left\| {\stackrel{\circ}{\ell}(x)|L_2^{(m)*} (0,1)} \right\|^2  =
( - 1)^{m + 1} \left[ {{{h^{2m} B_{2m} } \over {(2m)!}} + {{2h^{2m
+ 1} } \over {(2m)!}}\left\{ {\sum\limits_{\beta  = 0}^{t - 1}
{\left( {C_\beta h^{ - 1} \eta _\beta ^{2m}  - \beta ^{2m} }
\right) + } } \right.} \right.
$$
$$
+ \left. {\left. {\sum\limits_{k = 1}^{m - 1} {d_k \sum\limits_{i
= 0}^{2m} {{{( - 1)^i q_k^{t + i}  - q_k^{N - t + 1} } \over {(q_k
- 1)^{i+1}}}\Delta ^i t^{2m} } } } \right\}} \right],
$$
where $B_\alpha$ are Bernoulli numbers, $C_\beta$, $\beta  =
\overline {0,t - 1}$ are determined  from the  system (3.3), $q_k$
are the roots of the Euler-Frobenius  polynomial $E_{2m - 2}(q)$,
$\left| {q_k } \right| < 1$, $\eta _\beta  ,\,\,\,\beta  =
\overline {0,t - 1}$ are defined from (1.5), $\Delta ^i t^{2m} $
is finite-difference of order $i$ of $t^{2m}$, $d_k$ are
determined from (4.2). }

\textbf{Proof.} Taking into account (4.3), we reduce the
expression (1.7) to the form
    $$
\left\| {\ell(x) |L_2^{(m)*} (0,1)} \right\|^2  = ( - 1)^{m + 1}
\left[ {\sum\limits_{\beta  = 0}^N {C_\beta  \left( {f_m (x_\beta
) + P_{m - 1} (x_\beta )} \right) - {1 \over {(2m + 1)!}}} }
\right].\eqno (4.23)
$$
Applying the binomial formula to the equality (4.5), we obtain
    $$
f_m (x_\beta  ) = {{(1 - x_\beta  )^{2m} } \over {(2m)!}} -
\sum\limits_{i = 0}^{2m - 1} {{{( - 1)^{i + 1} x_\beta ^{2m - i -
1} } \over {2(i + 1)!(2m - i - 1)!}}.}\eqno (4.24)
$$

Putting (4.24) to the equality (4.23) and using equalities (4.17),
(1.9), consecutively we have
$$
\left\| {\ell(x) } \right\|^2  = ( - 1)^{m + 1} \left[
{\sum\limits_{\beta  = 0}^N {C_\beta  {{(1 - x_\beta  )^{2m} }
\over {(2m)!}}}  - {1 \over {(2m + 1)!}} - \sum\limits_{j = 1}^{m
- 1} {{{h^{2m - j} B_{2m - j} } \over {(j + 1)!(2m - j)!}} - } }
\right.
$$
$$
- \sum\limits_{j = 0}^{m - 1} {{{h^{2m - j} } \over {(j + 1)!(2m -
j - 1)!}}\left\{ {( - 1)^{2m - j - 1} \sum\limits_{\gamma  = 0}^{t
- 1} {\left( {C_\gamma  h^{ - 1} \eta _\gamma ^{2m - j - 1}  -
\gamma ^{2m - j - 1} } \right) + } } \right.}
$$
$$
+ \sum\limits_{k = 1}^{m - 1} {d_k \left. {\left. {\left(
{\sum\limits_{i = 1}^{2m - j - 1} {{{ - q_k  + ( - 1)^i q_k^{N +
i} } \over {(q_k  - 1)^{i + 1} }}\Delta ^i 0^j  + ( - 1)^{2m - j}
\sum\limits_{\gamma  = 1}^{t - 1} {\gamma ^{2m - j - 1} \left(
{q_k^\gamma   + q_k^{N - \gamma } } \right)} } } \right)}
\right\}} \right]}  =
$$
$$
 = ( - 1)^{m + 1} \left[ {Z_1  + Z_2 } \right],
 \eqno (4.25)
$$
where $Z_1  = \sum\limits_{\beta  = 0}^N {C_\beta  {{(1 - x_\beta
)^{2m} } \over {(2m)!}}}$ and $Z_2$ is the  remaining part in
square brackets of the equality (4.25).

Keeping in mind symmetry of the nodes (1.5), and doing similarly
calculations as in the proof of the lemma 3.3 when $\alpha = 2m$,
for $Z_1$ we have
$$
Z_1  = {1 \over {(2m + 1)!}} + {{h^{2m} B_{2m} } \over {(2m)!}} +
\sum\limits_{j = 1}^{2m - 2} {{{B_{2m - j} h^{2m - j} } \over {(j
+ 1)!(2m - j)!}} + }
$$
$$
+ \sum\limits_{j = 0}^{2m - 2} {{{h^{2m - j} } \over {(j + 1)!(2m
- j - 1)!}}\left[ {( - 1)^{2m - j - 1} \sum\limits_{\beta  = 0}^{t
- 1} {\left( {C_\beta  h^{ - 1} \eta _\beta ^{2m - j - 1}  - \beta
_{}^{2m - j - 1} } \right) + } } \right.}
$$
$$
+ \sum\limits_{k = 1}^{m - 1} {d_k \left. {\left\{ {\sum\limits_{i
= 1}^{2m - j - 1} {{{ - q_k^{N + i}  + ( - 1)^i q_k } \over {(1 -
q_k )^{i + 1} }}\Delta ^i 0^{2m - j - 1}  + ( - 1)^{2m - j}
\sum\limits_{\beta  = 0}^{t - 1} {\beta ^{2m - j - 1} \left(
{q_k^\beta   + q_k^{N - \beta } } \right)} } } \right\}} \right] +
}
$$
$$
+ {{2h^{2m + 1} } \over {(2m)!}}\left[ \sum\limits_{k = 1}^{m - 1}
d_k \left\{ {\sum\limits_{i = 1}^{2m} {{{ - q_k^{N + i}  + ( -
1)^i q_k } \over {(1 - q_k )^{i + 1} }}\Delta ^i 0^{2m}  -
\sum\limits_{\beta  = 0}^{t - 1} {\beta ^{2m} \left( {q_k^\beta +
q_k^{N - \beta } } \right)} } } \right\} + \right.
$$
$$
+\left. \sum\limits_{\beta  = 0}^{t - 1} {\left( {C_\beta  h^{ -
1} \eta _\beta ^{2m}  - \beta ^{2m} } \right)}   \right].\eqno
(4.26)
$$
Putting the equality (4.26) into (4.25) and taking into account
the equality (4.19), we have
$$
\left\| {\ell(x) } \right\|^2  = ( - 1)^{m + 1} \left[ {{{h^{2m}
B_{2m} } \over {(2m)!}} + {{2h^{2m + 1} } \over {(2m)!}}\left[
{\sum\limits_{k = 1}^{m - 1} {d_k \left\{ {\sum\limits_{i =
1}^{2m} {{{ - q_k^{N + i}  + ( - 1)^i q_k } \over {(1 - q_k )^{i +
1} }}\Delta ^i 0^{2m}  - } } \right.} } \right.} \right.
$$
$$
- \left. {\left. {\left. {\sum\limits_{\beta  = 1}^{t - 1} {\beta
^{2m} \left( {q_k^\beta   + q_k^{N - \beta } } \right)} } \right\}
+ \sum\limits_{\beta  = 0}^{t - 1} {\left( {C_\beta  h^{ - 1} \eta
_\beta ^{2m}  - \beta ^{2m} } \right)} } \right]} \right].
\eqno(4.27)
$$
Hence, using formulas (2.4), (2.5) and lemma 3.2 when $s = 0$,
after simplifications we will get the assertion of the theorem.

Theorem 4.2 is proved.

\textbf{Remark.} Should be noted, that when the nodes (1.5) are
equal spaced, i.e. when in (1.5) $\eta_0=0,\ \eta_1=1, ... ,
\eta_{t-1}=t-1$, from the theorems 4.1 and 4.2 we get the theorem
2.1 of [11] and the results of [31].

\section{Numerical results}

For simplicity in this section we investigate the optimal
quadrature formulas of the form (1.1) with the nodes
$$
x_0  = \eta _0 h,\,\,\,\,x_{N}  = 1 - \eta _0 h,\ \ \ 0 \le \eta
_0 <1,\ \  x_\beta   = h\beta,\ \ \beta=1,2,...,N-1. \eqno (5.1)
$$
This means we consider optimal quadrature formulas of the form
$$
\int\limits_0^1\varphi(x)dx\cong
C_0\varphi(\eta_0h)+\sum\limits_{\beta=1}^{N-1}C_{\beta}\varphi(h\beta)+
C_N\varphi(1-\eta_0h)\eqno (5.2)
$$
with the error functional
$$
\ell(x)=\varepsilon_{[0,1]}(x)-\left(C_0\delta(x-\eta_0h)+
\sum\limits_{\beta=1}^{N-1}C_{\beta}\delta(x-h\beta)+C_N\delta(x-(1-\eta_0h))\right).
\eqno (5.3)
$$
Then for this case from theorems 4.1 and 4.2 when $t=1$ we get
following

\textbf{Corollary 5.1} {\it The coefficients of optimal quadrature
formulas of the form (5.2) with the error functional (5.3) in the
Sobolev space $L_2^{(m)}(0,1)$ are expressed by formulas
$$
C_{\beta}=\left\{
\begin{array}{ll}
h\left(\frac{1}{2}+\sum\limits_{k=1}^{m-1}d_k\frac{q_k-q_k^N}{q_k-1}\right),\
& \beta=0,\ N,\\
h\left(1+\sum\limits_{k=1}^{m-1}d_k\left(q_k^{\beta}+q_k^{N-\beta}\right)\right),
&\beta=1,2,...,N-1, \\
\end{array}
\right. \eqno (5.4)
$$
where $d_k\ \ (k=\overline{1,m-1})$ satisfy following system of
$m-1$ linear equations
$$
\sum\limits_{k=1}^{m-1}d_k
\left(\sum\limits_{i=1}^j\frac{-q_k^2+(-1)^iq_k^{N-1+i}}{(q_k-1)^{i+1}}\Delta^i0^j+
\frac{(q_k-q_k^N)(1-\eta_0)^j}{q_k-1}\right)=
$$
$$
=\frac{1-B_{j+1}}{j+1}-\frac{(1-\eta_0)^j}{2}, \ \ \
j=1,2,...,m-1,\eqno (5.5)
$$
here $q_k$ are the roots of the Euler-Frobenius polynomial
$E_{2m-2}(q)$ of degree $2m-2$, $|q_k|<1$, $B_{j+1}$ are Bernoulli
numbers, $0\leq \eta_0<1$, $\Delta^i\gamma^j$ is finite difference
of order $i$ of $\gamma^j$,
$\Delta^i0^j=\Delta^i\gamma^j|_{\gamma=0}$.}

\textbf{Corollary 5.2.} {\it Square of norm of the error
functional (5.3) of optimal quadrature formulas of the form (5.2)
on the space $L_2^{(m)}(0,1)$ have following form
$$
\left\| {\stackrel{\circ}{\ell}(x)|L_2^{(m)*} (0,1)} \right\|^2  =
( - 1)^{m + 1} \left[ {{h^{2m} B_{2m} } \over {(2m)!}} + \right.
$$
$$
\left. +{{2h^{2m + 1} } \over {(2m)!}}\left(\frac{\eta_0^{2m}}{2}+
\sum\limits_{k=1}^{m-1}d_k\left[
\sum\limits_{i=1}^{2m}\frac{-q_k^{N+i}+(-1)^iq_k}{(1-q_k)^{i+1}}\Delta^i0^{2m}
+\frac{(q_k-q_k^N)\eta_0^{2m}}{q_k-1}\right]\right)\right],\eqno
(5.6)
$$
where $B_{2m}$ is Bernoulli number, $q_k$ are the roots of the
Euler-Frobenius  polynomial $E_{2m - 2}(q)$, $\left| {q_k }
\right| < 1$, $0\leq \eta_0<1$, $\Delta^i\gamma^{2m}$ is finite
difference of order $i$ of $\gamma^{2m}$,
$\Delta^i0^{2m}=\Delta^i\gamma^{2m}|_{\gamma=0}$. }

The aim of this section is by using the equalities (5.4), (5.5)
and by choosing $0\le \eta_0<1$ to investigate positiveness of
optimal coefficients.

As mentioned above, I.J.Schoenberg and S.Silliman in [25] showed,
that in the case of equal spaced nodes among optimal coefficients
it is appears negative coefficient starting from $m=7$. From (5.4)
and (5.5) in the case $\eta_0=0$ and $N\to \infty$ taking instead
of optimal coefficients $C_{\beta}$ the coefficients $NC_{\beta}$
we get the results of [25].

By choosing the value of $\eta_0=0,205$ and using (5.4) and (5.5),
we obtain optimal quadrature formulas of the form (5.2) with
positive coefficients in the cases $m=2,3,...,14$.

Below, in the first subsection,  we give the lists of the optimal
coefficients in tables. In the second subsection, we compared some
of the results of this paper with well known formulas.

\subsection{The coefficients of optimal quadrature formulas of the
form (5.2) in the case $\eta_0=0,205$}

Since $|q_k|<1$ from (5.4) it is easy to see, that the optimal
coefficients which located in the middle of the list of optimal
coefficients are close to $h$. Therefore, and also taking into
account symmetry of the optimal coefficients, we give the tables
of optimal coefficients which satisfy the condition
$|C_{\beta}-h|>10^{-8}$ and are located in the left boundary layer
in the list of optimal coefficients when $N=300$, $\eta_0=0,205$
for $m=2,3,...,14,15$.

\textbf{Table 1.} The optimal coefficients for $m=2$\\[1mm]
{\tiny
\begin{tabular}{|l|l|l|}
\hline
$C_0=0.00177612633884071008$& $C_3=0.00332336870971015344$ &$C_6=0.00333352503163128086$ \\
$C_1=0.00319454403039647954$& $C_4=0.00333600334618604447$ &$ C_7=0.00333328196792920787$\\
$C_2=0.00337052181497334175$ &$ C_5=0.00333261790554566866$& $ C_8=0.00333334709665188764$\\
\hline
\end{tabular}
 }\\[1mm]

\textbf{Table 2.} The optimal coefficients for $m=3$\\[1mm]
{\tiny
\begin{tabular}{|l|l|l|}
\hline
$C_0=0.00175247906611553142$& $C_3=0.00333236955552006215$ &$C_6=0.00333339747877195358$ \\
$C_1=0.00324264900163471085$& $C_4=0.00333368616893237390$ &$C_7=0.00333330571886293306$\\
$C_2=0.00333892060971646866$ &$C_5=0.00333318408917128015$& $C_8=0.00333334522322914572$\\
\hline
\end{tabular}
 }\\[1mm]

\textbf{Table 3.} The optimal coefficients for $m=4$\\[1mm]
{\tiny
\begin{tabular}{|l|l|l|}
\hline
$C_0=0.00174455012038852683$& $C_6 = 0.00332924765792690917$ &$C_{12} = 0.00333323709636443753$ \\
$C_1=0.00326634746379700186$& $C_7 = 0.00333552260766573748$ &$C_{13} = 0.00333338484710726915$\\
$C_2=0.00330774531081996009$ &$C_8 = 0.00333216117670071227$& $C_{14} = 0.00333330575901727716$\\
$C_3=0.00335699483364953318$ &$C_9 = 0.00333396080026924074$& $C_{15} = 0.00333334809332522744$\\
$C_4=0.00331942429918189420$ &$C_{10} = 0.00333299745834218190$& \\
$C_5 = 0.00334093116706873427$ &$C_{11} = 0.00333351312116039613$& \\
\hline
\end{tabular}
 }\\[1mm]

\textbf{Table 4.} The optimal coefficients for $m=5$\\[1mm]
{\tiny
\begin{tabular}{|l|l|l|}
\hline
$C_0 = 0.00173901538774817543 $& $  C_8 = 0.00332569334718594541 $ &$  C_{16} = 0.00333319051620654561 $ \\
$C_1 = 0.00328716890432268279  $& $  C_9 = 0.00333798181479271212 $ &$  C_{17} = 0.00333342016578285375 $ \\
$C_2 = 0.00326911260567167412  $& $  C_{10} = 0.00333050638482074723 $ &$  C_{18} = 0.00333328053942885252 $ \\
$  C_3 = 0.00340205381465105221  $& $  C_{11} = 0.00333505224863080752 $ &$  C_{19} = 0.00333336543188980042 $ \\
$  C_4 = 0.00328232898469587043  $& $  C_{12} = 0.00333228820948218315 $ &$  C_{20} = 0.00333331381749472885 $ \\
$  C_5 = 0.00336635018632553123  $& $  C_{13} = 0.00333396877152254017 $ &$  C_{21} = 0.00333334519891226771 $ \\
$  C_6 = 0.00331284808226362786  $& $  C_{14} = 0.00333294698744021598 $ &$ $ \\
$  C_7 = 0.00334587152301716133  $& $  C_{15} = 0.00333356823085626853 $ &$ $ \\
\hline
\end{tabular}
 }\\[0.5mm]

\textbf{Table 5.} The optimal coefficients for $m=6$\\[1mm]
{\tiny
\begin{tabular}{|l|l|l|l|}
\hline
$C_0 = 0.00173498886058381026 $&$C_7 = 0.00337747643535115270 $ &$C_{15} = 0.00333497426498772165  $ &$C_{21} = 0.00333347053324047546$ \\
$C_1 = 0.00330520593286011352 $&$C_8 = 0.00330385537900635393 $&$C_{14} = 0.00333085188184676641$  &$C_{22} = 0.00333324260768654442$ \\
$C_2 = 0.00322620390097277555 $&$C_9 = 0.00335290464626403816 $ &$C_{16} = 0.00333224823220328272 $ &$C_{23} = 0.00333339332712609729$ \\
$C_3 = 0.00346783504823977116 $& $C_{10} = 0.00332037008745176751  $ &$C_{17} = 0.00333405087626146100  $ &$C_{24} = 0.00333329366147357165$ \\
$C_4 = 0.00321158077766250064 $& $C_{11} = 0.00334191131490927107  $ & $C_{18} = 0.00333285884589705156  $ & $C_{25} = 0.00333335956698815529$ \\
$C_5 = 0.00342625617820782061 $& $C_{12} =0.003327659419014779572  $ & $C_{19} = 0.00333364709595060379$  & $C_{26} = 0.00333331598590751689$ \\
$C_6 = 0.00326815998330420462 $& $C_{13} = 0.00333708573208910388  $& $C_{20} = 0.00333312585271306499$  & $C_{27} = 0.00333334480459741380$ \\
\hline
\end{tabular}
 }\\[1mm]

\textbf{Table 6.} The optimal coefficients for $m=7$\\[1mm]
{\tiny
\begin{tabular}{|l|l|l|l|}
\hline
$C_0 = 0.00173205877707359127  $&$C_9 = 0.00339643943424144824  $  & $C_{18} = 0.00333069296984364163   $ &$C_{27} = 0.00333344250631211799   $\\
$C_1 = 0.00332034204314634415  $& $C_{10} = 0.00328876219893682112  $ & $C_{19} = 0.00333518660336370406   $ &$C_{28} = 0.00333325670544635996   $\\
$C_2 = 0.00318230520695443787  $&$C_{11} = 0.00336471020635444449  $ & $C_{20} = 0.00333203252906371234   $ &$C_{29} = 0.00333338711800696554   $ \\
$C_3 = 0.00355132503340265203  $& $C_{12} = 0.00331127918647367579   $&  $C_{21} = 0.00333424636193460625   $&$C_{30} = 0.00333329558217999618   $\\
$C_4 = 0.00309977448873935734  $& $C_{13} = 0.00334882331050388661   $ &  $C_{22} = 0.00333269248328614613   $ &$C_{31} = 0.00333335983065088640   $ \\
$C_5 = 0.00354150500848728542  $& $C_{14} = 0.00332245757624453130   $ & $C_{23} = 0.00333378314247101608   $  &$C_{32} = 0.00333331473501844551   $ \\
$C_6 = 0.00316818905922230260  $&$C_{15} = 0.00334096810776062514   $ &$C_{24} = 0.00333301761482745647   $  &$C_{33} = 0.00333334638738364948   $ \\
$C_7 = 0.00345627307134070750  $& $C_{16} = 0.00332797414827915069   $&$C_{25} = 0.00333355493435702833   $&        \\
$C_8 = 0.00324460285490503001  $&$C_{17} = 0.00333709504135639172   $ &$C_{26} = 0.00333317779284220016   $ &       \\
\hline
\end{tabular}
 }\\[1mm]

\textbf{Table 7.} The optimal coefficients for $m=8$\\[1mm]
{\tiny
\begin{tabular}{|l|l|l|l|}
\hline
$C_0 = 0.00172995290066120350  $&$C_{11} = 0.00342842342134998057 $  & $C_{22} = 0.00333014109574784269 $ &$C_{33} = 0.00333343948477634317 $\\
$C_1 = 0.00333261654096956311  $&$C_{12} = 0.00326323297211348103  $  & $C_{23} = 0.00333567603784529139 $ &$C_{34} = 0.00333325543170061418 $\\
$C_2 = 0.00314030615358630034  $&$C_{13} = 0.00338490036159433393  $  & $C_{24} = 0.00333161408331454461 $ &$C_{35} = 0.00333339050320497157 $\\
$C_3 = 0.00364680181687277409  $&$C_{14} = 0.00329544249146043423  $  & $C_{25} = 0.00333459504509005457 $ &$C_{36} = 0.00333329137793257419 $\\
$C_4 = 0.00294599819613147818  $&$C_{15} = 0.00336115859898192848  $  & $C_{26} = 0.00333240739717116319 $ &$C_{37} = 0.00333336412325119880 $\\
$C_5 = 0.00373067122924529534  $&$C_{16} = 0.00331290610817136431  $  & $C_{27} = 0.00333401285268273812 $ &$C_{38} = 0.00333331073745712191 $\\
$C_6 = 0.00297727240591131934  $&$C_{17} = 0.00334832702306404860  $  & $C_{28} = 0.00333283465264497077 $ &$C_{39} = 0.00333334991582712020 $\\
$C_7 = 0.00362534555471145307  $&$C_{18} = 0.00332232883088807872  $  & $C_{29} = 0.00333369930144181951 $ &$C_{40} = 0.00333332116389597219 $\\
$C_8 = 0.00310592534617643848  $&$C_{19} = 0.00334140963871753226  $  & $C_{30} = 0.00333306475936511862 $ &$ $\\
$C_9 = 0.00350552648624077583  $&$C_{20} = 0.00332740619898615369  $  & $C_{31} = 0.00333353043240642612 $ &$ $\\
$C_{10} = 0.00320487233533179992  $&$C_21 = 0.00333768315455876029  $  & $C_{32} = 0.00333318868772800817 $ &$ $\\
\hline
\end{tabular}
 }\\[1mm]

\textbf{Table 8.} The optimal coefficients for $m=9$\\[1mm]
{\tiny
\begin{tabular}{|l|l|l|l|}
\hline
$C_0 = 0.00172847784910360447  $&$C_{12} = 0.00314035977572450425  $  & $C_{24} = 0.00332610031017101433  $ &$C_{36} = 0.00333306608365147515 $\\
$C_1 = 0.00334216427951079832  $&$C_{13} = 0.00348084679101216452  $  & $C_{25} = 0.00333882817755420089  $ &$C_{37} = 0.00333353635843483279 $\\
$C_2 = 0.00310265417511838644  $&$C_{14} = 0.00322087398274449081  $  & $C_{26} = 0.00332915897604222655  $ &$C_{38} = 0.00333317909856898933 $\\
$C_3 = 0.00374645947596537845  $&$C_{15} = 0.00341893779189990479  $  & $C_{27} = 0.00333650452957059022  $ &$C_{39} = 0.00333345050289666260 $\\
$C_4 = 0.00275792168171570242  $&$C_{16} = 0.00326822729185132696  $  & $C_{28} = 0.00333092422572131367  $ &$C_{40} = 0.00333324432158157594 $\\
$C_5 = 0.00400138156396844176  $&$C_{17} = 0.00338282511219261974  $  & $C_{29} = 0.00333516349327706089  $ &$C_{41} = 0.00333340095406794224 $\\
$C_6 = 0.00266165221372250548  $&$C_{18} = 0.00329572154562808534  $  & $C_{30} = 0.00333194299096395421  $ &$C_{42} = 0.00333328196299562843 $\\
$C_7 = 0.00394079887384133675  $&$C_{19} = 0.00336191228831355167  $  & $C_{31} = 0.00333438955334958306  $ &$C_{43} = 0.00333337235851711259 $\\
$C_8 = 0.00282156182791894007  $&$C_{20} = 0.00331161983320458289  $  & $C_{32} = 0.00333253094060154966  $ &$C_{44} = 0.00333330368655589561 $\\
$C_9 = 0.00374603753952685103  $&$C_{21} = 0.00334982981298158293  $  & $C_{33} = 0.00333394289767582404  $ &$C_{45} = 0.00333335585549280241 $\\
$C_{10} = 0.00300894501757551301  $&$C_{22} = 0.00332080076136125273  $  & $C_{34} = 0.00333287025750894230  $ &$C_{46} = 0.00333331622362661400 $\\
$C_{11} = 0.00358457992357594231  $&$C_{23} = 0.00334285432190737266  $  & $C_{35} = 0.00333368512427696033  $ &$C_{47} = 0.00333334633129049714 $\\
\hline
\end{tabular}
 }\\[1mm]

\textbf{Table 9.} The optimal coefficients for $m=10$\\[1mm]
{\tiny
\begin{tabular}{|l|l|l|l|}
\hline
$C_0 = 0.00172749278862645234  $&$C_{14} = 0.00304088010987428129   $  & $C_{28} = 0.00332404459372417715   $ &$C_{42} = 0.00333304181986959795  $\\
$C_1 = 0.00334915954552527705  $&$C_{15} = 0.00356281566915249473   $  & $C_{29} = 0.00334058737008624209   $ &$C_{43} = 0.00333356098973720558  $\\
$C_2 = 0.00307141154907661477  $&$C_{16} = 0.00315360445472180508   $  & $C_{30} = 0.00332766830567799599   $ &$C_{44} = 0.00333315554587373269  $\\
$C_3 = 0.00384087716305041600  $&$C_{17} = 0.00347393491400228088   $  & $C_{31} = 0.00333775742286366819   $ &$C_{45} = 0.00333347217581801177  $\\
$C_4 = 0.00255327063613243166  $&$C_{18} = 0.00322341661691745051   $  & $C_{32} = 0.00332987835336878476   $ &$C_{46} = 0.00333322490478845213  $\\
$C_5 = 0.00434013033335869961  $&$C_{19} = 0.00341922625493328018   $  & $C_{33} = 0.00333603148905263082   $ &$C_{47} = 0.00333341801022001505  $\\
$C_6 = 0.00221013903097713541  $&$C_{20} = 0.00326623022080639272   $  & $C_{34} = 0.00333122621755866373   $ &$C_{48} = 0.00333326720521970105  $\\
$C_7 = 0.00444946400666106578  $&$C_{21} = 0.00338574919280914013   $  & $C_{35} = 0.00333497887821597101   $ &$C_{49} = 0.00333338497584858724  $\\
$C_8 = 0.00231492992409806705  $&$C_{22} = 0.00329239374555605829   $  & $C_{36} = 0.00333204825076180440   $ &$C_{50} = 0.00333329300329463880  $\\
$C_9 = 0.00420985817056583358  $&$C_{23} = 0.00336530759722583348   $  & $C_{37} = 0.00333433691406739092   $ &$C_{51} = 0.00333336482893367433  $\\
$C_{10} = 0.00260680747387077026  $&$C_{24} = 0.00330836190455046297   $  & $C_{38} = 0.00333254959049473769   $ &$C_{52} = 0.00333330873695621772  $\\
$C_{11} = 0.00392158564879464832  $&$C_{25} = 0.00335283526564426801   $  & $C_{39} = 0.00333394539453196409   $ &$C_{53} = 0.00333335254178659839  $\\
$C_{12} = 0.00286380539858292753  $&$C_{26} = 0.00331810309402904888   $  & $C_{40} = 0.00333285534630910409   $ &$C_{54} = 0.00333331833256004469  $\\
$C_{13} = 0.00370486410231191268  $&$C_{27} = 0.00334522746373694056   $  & $C_{41} = 0.00333370661560394268   $ &$C_{55} = 0.00333334504813390320  $\\
\hline
\end{tabular}
 }\\[1mm]

\textbf{Table 10.} The optimal coefficients for $m=11$\\[1mm]
{\tiny
\begin{tabular}{|l|l|l|l|}
\hline
$C_0 = 0.00172689334490066565  $&$  C_{16} = 0.00290248747287229047    $  & $  C_{32} = 0.00332138043174862451    $ &$  C_{48} = 0.00333300511542151927   $\\
$C_1 = 0.00335378561374918185  $&$ C_{17} = 0.00367878243460793137    $  & $  C_{33} = 0.00334288089094969269    $ &$  C_{49} = 0.00333359550162011932   $\\
$C_2 = 0.00304832508611045170  $&$  C_{18} = 0.00305674037361827340    $  & $  C_{34} = 0.00332570708681864833    $ &$  C_{50} = 0.00333312392303583896   $\\
$C_3 = 0.00391935697738649923  $&$  C_{19} = 0.00355459988011267817    $  & $  C_{35} = 0.00333942490198512218    $ &$  C_{51} = 0.00333350060250695466   $\\
$  C_4 = 0.00236091789799528996  $&$  C_{20} = 0.00315642482727045543    $  & $  C_{36} = 0.00332846761136009168    $ &$  C_{52} = 0.00333319972492547481   $\\
$  C_5 = 0.00470103221229688833  $&$  C_{21} = 0.00347472666712902671    $  & $  C_{37} = 0.00333721989298164747    $ &$  C_{53} = 0.00333344005477077078   $\\
$  C_6 = 0.00166606164689116412  $&$  C_{22} = 0.00322035045788260487    $  & $  C_{38} = 0.00333022889296381827    $ &$  C_{54} = 0.00333324808820998967   $\\
$  C_7 = 0.00513715438246413703  $&$  C_{23} = 0.00342360164155868860    $  & $  C_{39} = 0.00333581304546978635    $ &$  C_{55} = 0.00333340142397405117   $\\
$  C_8 = 0.00155687727226540244  $&$  C_{24} = 0.00326121936040494949    $  & $  C_{40} = 0.00333135263111543830    $ &$  C_{56} = 0.00333327894505666322   $\\
$  C_9 = 0.00496522549139354527  $&$  C_{25} = 0.00339094085206410810    $  & $  C_{41} = 0.00333491544485158395    $ &$  C_{57} = 0.00333337677667015449   $\\
$  C_{10} = 0.00190545757664787945  $&$   C_{26} = 0.00328731578129565660   $  & $  C_{42} = 0.00333206960132156587    $ &$  C_{58} = 0.00333329863240825021   $\\
$  C_{11} = 0.00454262920324724638  $&$  C_{27} = 0.00337009186182822834    $  & $  C_{43} = 0.00333434275558700776    $ &$  C_{59} = 0.00333336105114438844   $\\
$  C_{12} = 0.00233068337371350425  $&$  C_{28} = 0.00330397127017939317    $  & $  C_{44} = 0.00333252704428838619    $ &$  C_{60} = 0.00333331119337459038   $\\
$  C_{13} = 0.00415339159536806206  $&$  C_{29} = 0.00335678702055923429    $  & $  C_{45} = 0.00333397736710399180    $ &$  C_{61} = 0.00333335101791066748   $\\
$  C_{14} = 0.00266841337297767427  $&$  C_{30} = 0.00331459921621225224    $  & $  C_{46} = 0.00333281890305696659    $ &$   C_{62} = 0.00333331920754931007  $\\
$  C_{15} = 0.00386950335147916675  $&$  C_{31} = 0.00334829754279374992    $  & $  C_{47} = 0.00333374424114699843    $ &$  C_{63} = 0.00333334461648534207   $\\
\hline
\end{tabular}
 }\\[1mm]

\textbf{Table 11.} The optimal coefficients for $m=12$\\[1mm]
{\tiny
\begin{tabular}{|l|l|l|l|}
\hline
$  C_0 = 0.00172660093000204434   $&$  C_{18} = 0.00275760059774337822       $  & $  C_{36} = 0.00331906522131345626     $ &$  C_{54} = 0.00333298282252256406    $\\
$  C_1 = 0.00335622016635424541   $&$  C_{19} = 0.00380333655344864095       $  & $  C_{37} = 0.00334494624845792310     $ &$  C_{55} = 0.00333361861605357889    $\\
$  C_2 = 0.00303488707201474333   $&$  C_{20} = 0.00295003782640532812       $  & $  C_{38} = 0.00332388150649666917     $ &$  C_{56} = 0.00333310114012289132    $\\
$  C_3 = 0.00397016939076339731   $&$  C_{21} = 0.00364570570477646406       $  & $  C_{39} = 0.00334102623365402826     $ &$  C_{57} = 0.00333352231669313408    $\\
$  C_4 = 0.00222172391357403828   $&$  C_{22} = 0.00307887414398914889       $  & $  C_{40} = 0.00332707203566076212     $ &$  C_{58} = 0.00333317951871399489    $\\
$  C_5 = 0.00499378828442308939   $&$  C_{23} = 0.00354055563717261057       $  & $  C_{41} = 0.00333842944024015643     $ &$  C_{59} = 0.00333345852391701871    $\\
$  C_6 = 0.00117134786383522544   $&$  C_{24} = 0.00316461152396595462       $  & $  C_{42} = 0.00332918558244974745     $ &$  C_{60} = 0.00333323144001324420    $\\
$  C_7 = 0.00583516570203082313   $&$  C_{25} = 0.00347069038125884784       $  & $  C_{43} = 0.00333670921153828293     $ &$  C_{61} = 0.00333341626487956293    $\\
$  C_8 = 0.00070537581513706783   $&$  C_{26} = 0.00322151976711594156       $  & $  C_{44} = 0.00333058568684464716     $ &$  C_{62} = 0.00333326583488148334    $\\
$  C_9 = 0.00589276164637635086   $&$  C_{27} = 0.00342434860282943662       $  & $  C_{45} = 0.00333556965838687970     $ &$  C_{63} = 0.00333338827070287737    $\\
$  C_{10} = 0.00097691145372704574   $&$  C_{28} = 0.00325925036543305965       $  & $  C_{46} = 0.00333151317572930858     $ &$  C_{64} = 0.00333328861949931687    $\\
$  C_{11} = 0.00541764827809579119   $&$  C_{29} = 0.00339363262285636622       $  & $  C_{47} = 0.00333481476987879354     $ &$  C_{65} = 0.00333336972617429080    $\\
$  C_{12} = 0.00154000801804920969   $&$  C_{30} = 0.00328425393283165749       $  & $  C_{48} = 0.00333212758377888206     $ &$  C_{66} = 0.00333330371299542040    $\\
$  C_{13} = 0.00484769113306215419   $&$  C_{31} = 0.00337328012014849898       $  & $  C_{49} = 0.00333431469968656027     $ &$  C_{67} = 0.00333335744149293259    $\\
$  C_{14} = 0.00207043279270397406   $&$  C_{32} = 0.00330081999773997306       $  & $  C_{50} = 0.00333253459373417958     $ &$  C_{68} = 0.00333331371156666169    $\\
$  C_{15} = 0.00437783417724906036   $&$  C_{33} = 0.00335979634546885836       $  & $  C_{51} = 0.00333398343198778386     $ &$  C_{69} = 0.00333334930359957163    $\\
$  C_{16} = 0.00247418267818456252   $&$  C_{34} = 0.00331179481824981668       $  & $  C_{52} = 0.00333280421438103343     $ &$  C_{70} = 0.00333332033504383217    $\\
$  C_{17} = 0.00403747731479976406   $&$  C_{35} = 0.00335086371649026079       $  & $  C_{53} = 0.00333376398622498520     $ &$  C_{71} = 0.00333334391271425501    $\\
\hline
\end{tabular}
 }\\[1mm]

\textbf{Table 12.} The optimal coefficients for $m=13$\\[1mm]
{\tiny
\begin{tabular}{|l|l|l|l|}
\hline
$ C_0 = 0.00172655548945567044  $&$  C_{20} = 0.00277237192137831871   $  & $  C_{40} = 0.00332071149652304087    $ &$  C_{60} = 0.00333305114925564172      $\\
$  C_1 = 0.00335662901376447258  $&$  C_{21} = 0.00379816301931711370   $  & $  C_{41} = 0.00334377075895046149    $ &$  C_{61} = 0.00333356668071934869      $\\
$  C_2 = 0.00303238425807657517  $&$  C_{22} = 0.00294841399870425801   $  & $  C_{42} = 0.00332470227406376349    $ &$  C_{62} = 0.00333314037062970202      $\\
$  C_3 = 0.00398073517633032120  $&$  C_{23} = 0.00365193776947454155   $  & $  C_{43} = 0.00334047064592444503    $ &$  C_{63} = 0.00333349290060179606      $\\
$  C_4 = 0.00218920482037733280  $&$  C_{24} = 0.00306969888569782702   $  & $  C_{44} = 0.00332743125080413949    $ &$  C_{64} = 0.00333320138185937509      $\\
$  C_5 = 0.00507104524959419857  $&$  C_{25} = 0.00355143697213385530   $  & $  C_{45} = 0.00333821396235558649    $ &$  C_{65} = 0.00333344244838961507      $\\
$  C_6 = 0.00102336156051060256  $&$  C_{26} = 0.00315292244059992812   $  & $  C_{46} = 0.00332929737855920021    $ &$  C_{66} = 0.00333324310246929006      $\\
$  C_7 = 0.00607199798241084696  $&$  C_{27} = 0.00348255132276115587   $  & $  C_{47} = 0.00333667079868275828    $ &$  C_{67} = 0.00333340794823161848      $\\
$  C_8 = 0.00037867881870999653  $&$  C_{28} = 0.00320992302093593156   $  & $  C_{48} = 0.00333057347222005968    $ &$  C_{68} = 0.00333327163179578226      $\\
$  C_9 = 0.00629239692582258769  $&$  C_{29} = 0.00343539498005851424   $  & $  C_{49} = 0.00333561555424564661    $ &$   C_{69} = 0.00333338435638400590     $\\
$  C_{10} = 0.00053211065682692867  $&$  C_{30} = 0.00324892979162245131   $  & $  C_{50} = 0.00333144608900867865    $ &$  C_{70} = 0.00333329114067812168      $\\
$  C_{11} = 0.00587832107921238828  $&$  C_{31} = 0.00340313244846696006   $  & $  C_{51} = 0.00333489395837593862    $ &$  C_{71} = 0.00333336822384125740      $\\
$  C_{12} = 0.00108774569344314697  $&$  C_{32} = 0.00327561242680296621   $  & $  C_{52} = 0.00333204280063615644    $ &$  C_{72} = 0.00333330448121352504      $\\
\hline
\end{tabular} }
\newpage
\textbf{Table 12.} Continuation\\[1mm]
{\tiny
\begin{tabular}{|l|l|l|l|}
\hline
$  C_{13} = 0.00527471154107237576  $&$  C_{33} = 0.00338106562394459128   $  & $  C_{53} = 0.00333440051769216630    $ &$  C_{73} = 0.00333335719210945782      $\\
$  C_{14} = 0.00167844383373434332  $&$  C_{34} = 0.00329386137938002601   $  & $  C_{54} = 0.00333245084313764205    $ &$  C_{74} = 0.00333331360371894311      $\\
$  C_{15} = 0.00473043759244898944  $&$  C_{35} = 0.00336597430528210186   $  & $  C_{55} = 0.00333406309378899346    $ &$ C_{75} = 0.00333334964840674512      $\\
$  C_{16} = 0.00216162520663729902  $&$  C_{36} = 0.00330634125957574669   $  & $  C_{56} = 0.00333272987017358315    $ &$  C_{76} = 0.00333331984185726252      $\\
$   C_{17} = 0.00431159258197410437 $&$  C_{37} = 0.00335565407117066211   $  & $  C_{57} = 0.00333383235706368672    $ &$  C_{77} = 0.00333334448988298376      $\\
$  C_{18} = 0.00251908554798301532  $&$  C_{38} = 0.00331487552014629517   $  & $  C_{58} = 0.00333292067403627449    $ &$      $\\
$  C_{19} = 0.00400965449172913550  $&$  C_{39} = 0.00334859674230775176   $  & $  C_{59} = 0.00333367457501187563    $ &$      $\\
\hline
\end{tabular}
 }\\[1mm]

\textbf{Table 13.} The optimal coefficients for $m=14$\\[1mm]
{\tiny
\begin{tabular}{|l|l|l|l|}
\hline
$  C_0 = 0.00172671049778565797  $&$  C_{23} = 0.00316025452220934915    $  & $  C_{46} = 0.00333638193499156475   $ &$  C_{69} = 0.00333328063425493357     $\\
$  C_1 = 0.00335516395507917580  $&$  C_{24} = 0.00347924238639198785    $  & $  C_{47} = 0.00333077781143822645   $ &$  C_{70} = 0.00333337750881443832     $\\
$  C_2 = 0.00304193627584372803  $&$  C_{25} = 0.00321053936117372331    $  & $  C_{48} = 0.00333547552503919473   $ &$  C_{71} = 0.00333329630283508029     $\\
$  C_3 = 0.00393776327125170017  $&$  C_{26} = 0.00343655097788938971    $  & $  C_{49} = 0.00333153762035049978   $ &$  C_{72} = 0.00333336437448481488     $\\
$   C_4 = 0.00233007712879415692 $&$  C_{27} = 0.00324664320823834275    $  & $  C_{50} = 0.00333483860682981228   $ &$  C_{73} = 0.00333330731280618138     $\\
$  C_5 = 0.00471527133757507641  $&$   C_{28} = 0.00340610011591843955   $  & $  C_{51} = 0.00333207152369661671   $ &$  C_{74} = 0.00333335514527654166     $\\
$  C_6 = 0.00174497854628782157  $&$  C_{29} = 0.00327227842278880512    $  & $  C_{52} = 0.00333439105698931688   $ &$  C_{75} = 0.00333331504927383940     $\\
$  C_7 = 0.00485675294324682823  $&$  C_{30} = 0.00338454689233366739    $  & $  C_{53} = 0.00333244668671512832   $ &$  C_{76} = 0.00333334866011226255     $\\
$  C_8 = 0.00212747598467839435  $&$  C_{31} = 0.00329038332124877627    $  & $  C_{54} = 0.00333407657300585164   $ &$  C_{77} = 0.00333332048552169609     $\\
$  C_9 = 0.00408453051443764890  $&$  C_{32} = 0.00336934817236508896    $  & $  C_{55} = 0.00333271030587360587   $ &$  C_{78} = 0.00333334410312812365     $\\
$  C_{10} = 0.00303849317411607795  $&$  C_{33} = 0.00330313676378294222    $  & $  C_{56} = 0.00333385559179154898   $ &$  C_{79} = 0.00333332430545503779     $\\
$  C_{11} = 0.00326150825980108299  $&$  C_{34} = 0.00335864987188439309    $  & $  C_{57} = 0.00333289554541894104   $ &$  C_{80} = 0.00333334090103430680     $\\
$  C_{12} = 0.00365086987677036012  $&$  C_{35} = 0.00331210917522999367    $  & $  C_{58} = 0.00333370031303950707   $ &$  C_{81} = 0.00333332698963926159     $\\
$  C_{13} = 0.00288027997265296082  $&$  C_{36} = 0.00335112604866268579    $  & $  C_{59} = 0.00333302570925529457   $ &$  C_{82} = 0.00333333865099250313     $\\
$  C_{14} = 0.00383972630105596425  $&$  C_{37} = 0.00331841762354822903    $  & $  C_{60} = 0.00333359120201437040   $ &$  C_{83} = 0.00333332887575723863     $\\
$  C_{15} = 0.00282694633373975373  $&$  C_{38} = 0.00334583703432045650    $  & $  C_{61} = 0.00333311717256726671   $ &$  C_{84} = 0.00333333706993697896     $\\
$  C_{16} = 0.00380915983064439615  $&$  C_{39} = 0.00332285171576454862    $  & $  C_{62} = 0.00333351453205997142   $ &$  C_{85} = 0.00333333020109139106     $\\
$  C_{17} = 0.00290295426940039930  $&$  C_{40} = 0.00334211980668418472    $  & $  C_{63} = 0.00333318144185430803   $ &$  C_{86} = 0.00333333595896358890     $\\
$  C_{18} = 0.00371316819510743415  $&$  C_{41} = 0.00332596789711823063    $  & $  C_{64} = 0.00333346065774556324   $ &$  C_{87} = 0.00333333113237492521     $\\
$  C_{19} = 0.00300351110208944757  $&$   C_{42} = 0.00333950753291683119   $  & $  C_{65} = 0.00333322660248814852   $ &$  C_{88} = 0.00333333517830670844     $\\
$  C_{20} = 0.00361660703987917620  $&$  C_{43} = 0.00332815772174759486    $  & $  C_{66} = 0.00333342280143380223   $ &$     $\\
$  C_{21} = 0.00309185019337987811  $&$  C_{44} = 0.00333767185593130795    $  & $  C_{67} = 0.00333325833588546756   $ &$     $\\
$  C_{22} = 0.00353813630227038258  $&$  C_{45} = 0.00332969651611532286    $  & $  C_{68} = 0.00333339620062519364   $ &$     $\\
\hline
\end{tabular} } \\[1mm]

\textbf{Table 14.} The optimal coefficients for $m=15$\\[1mm]
{\tiny
\begin{tabular}{|l|l|l|l|}
\hline
$  C_0 = 0.00172702945455083193  $&$  C_{26} = 0.00591821410669587006    $  & $  C_{52} = 0.00336928787611088775   $ &$  C_{78} = 0.00333383067161036060       $\\
$  C_1 = 0.00335196263611974402  $&$  C_{27} = 0.00113739105014957971    $  & $  C_{53} = 0.00330283683648456223   $ &$  C_{79} = 0.00333291149314649690       $\\
$  C_2 = 0.00306452563504523993  $&$  C_{28} = 0.00519802682335201819    $  & $  C_{54} = 0.00335920033583227050   $ &$  C_{80} = 0.00333369113636488924       $\\
$  C_3 = 0.00382735744574236094  $&$  C_{29} = 0.00175042759986189166    $  & $  C_{55} = 0.00331139305001559083   $ &$  C_{81} = 0.00333302984634137357       $\\
$  C_4 = 0.00272471104016856072  $&$  C_{30} = 0.00467672632259750788    $  & $  C_{56} = 0.00335194298912504867   $ &$  C_{82} = 0.00333359074969419812       $\\
$  C_5 = 0.00362523740203732478  $&$  C_{31} = 0.00219339925074939015    $  & $  C_{57} = 0.00331754870221566789   $ &$  C_{83} = 0.00333311499388354665       $\\
$  C_6 = 0.00416828584678143709  $&$  C_{32} = 0.00430050875356561753    $  & $  C_{58} = 0.00334672179022929618   $ &$  C_{84} = 0.00333351852791527129       $\\
$  C_7 = 0.00038025393759564147  $&$  C_{33} = 0.00251280401810681106    $  & $  C_{59} = 0.00332197730137501015   $ &$  C_{85} = 0.00333317625208569821       $\\
$  C_8 = 0.00918605541537786417  $&$  C_{34} = 0.00402940935773975411    $  & $  C_{60} = 0.00334296547103844503   $ &$  C_{86} = 0.00333346656897177222       $\\
$  C_9 = -0.0056482216521837456  $&$  C_{35} = 0.00274285983082156931    $  & $  C_{61} = 0.00332516339541746467   $ &$  C_{87} = 0.00333322032343801599       $\\
$  C_{10} = 0.01503979792439844429  $&$  C_{36} = 0.00383420982729031424    $  & $  C_{62} = 0.00334026303961781421   $ &$  C_{88} = 0.00333342918784036036       $\\
$  C_{11} = -0.0102513216994999018  $&$  C_{37} = 0.00290846808379631022    $  & $  C_{63} = 0.00332745558621977111   $ &$  C_{89} = 0.00333325202995283771       $\\
$  C_{12} = 0.01780870154874396918  $&$  C_{38} = 0.00369371677971561177    $  & $  C_{64} = 0.00333831881318149643   $ &$  C_{90} = 0.00333340229451091797       $\\
$  C_{13} = -0.0111516599198095529  $&$  C_{39} = 0.00302764878602438014    $  & $  C_{65} = 0.00332910467075450735   $ &$  C_{91} = 0.00333327484075914716       $\\
$  C_{14} = 0.01717584313564491950  $&$  C_{40} = 0.00359261903731864099    $  & $  C_{66} = 0.00333692006675210665   $ &$  C_{92} = 0.00333338294648316812       $\\
$  C_{15} = -0.0094596822947696162  $&$  C_{41} = 0.00311340503994861123    $  & $  C_{67} = 0.00333029108147133945   $ &$  C_{93} = 0.00333329125167301758       $\\
$  C_{16} = 0.01487372890505018792  $&$  C_{42} = 0.00351987755666928677    $  & $  C_{68} = 0.00333591375827837438   $ &$  C_{94} = 0.00333336902681665769       $\\
$  C_{17} = -0.0068972585187668098  $&$  C_{43} = 0.00317510611508734637    $  & $  C_{69} = 0.00333114462799414914   $ &$  C_{95} = 0.00333330305827423077       $\\
\hline
\end{tabular}
}
\newpage
\textbf{Table 14.} Continuation\\[1mm]
{\tiny
\begin{tabular}{|l|l|l|l|}
\hline
$  C_{18} = 0.01228998089227584958  $&$  C_{44} = 0.00346754171514042525    $  & $  C_{70} = 0.00333518978378703653   $ &$  C_{96} = 0.00333335901250834411       $\\
$  C_{19} = -0.0044379368869592209  $&$  C_{45} = 0.00321949791378515819    $  & $  C_{71} = 0.00333175870004468364   $ &$  C_{97} = 0.00333331155236741126       $\\
$  C_{20} = 0.01003287196894858974  $&$  C_{46} = 0.00342988831851899784    $  & $  C_{72} = 0.00333466893051466329   $ &$  C_{98} = 0.00333335180785510696       $\\
$  C_{21} = -0.0024157213542609094  $&$  C_{47} = 0.00325143564488660772    $  & $  C_{73} = 0.00333220048555005048   $ &$  C_{99} = 0.00333331766332322668       $\\
$  C_{22} = 0.00825049420457872496  $&$  C_{48} = 0.00340279870089829129    $  & $  C_{74} = 0.00333429420994935355   $ &$  C_{100} = 0.0033333466245687503       $\\
$  C_{23} = -0.0008623772730946671  $&$  C_{49} = 0.00327441304627985970    $  & $  C_{75} = 0.00333251832192158444   $ &$  C_{101} = 0.0033333220597643263       $\\
$  C_{24} = 0.00690738675157999725  $&$  C_{50}= 0.003383309304493717340    $  & $  C_{76} = 0.00333402462250633817   $ &$       $\\
$  C_{25} = 0.00029252137237138985  $&$  C_{51} = 0.00329094390700659850    $  & $  C_{77} = 0.00333274698483095824   $ &$       $\\
\hline
\end{tabular} } \\[1mm]

The tables 1-13 confirm positiveness of the optimal coefficients
of optimal quadrature formulas of the form (5.2) with the nodes
(5.1) when $\eta_0=0,205$.

\subsection{Comparison some of the results with well known formulas}

In this section we compare some of the results of this work with
the example (e) of [12]. From the example (e) of [12] when
$[a,b]=[0,1]$ and for the equidistant nodes
$$
x_{\beta}=h\beta,\ \ \beta=0,1,...,N,\ \ h=\frac{1}{N}
$$
of the interval $[0,1]$ we obtain following particular formulas
which were already obtained by Sard in [14, 18].

If $N=2$ the optimal quadrature formula of the type (1.1) is
$$
\int\limits_0^1\varphi(x)dx=\frac{1}{16}\left[3\varphi(0)+10\varphi(1/2)+3\varphi(1)\right]+
R[\varphi] \eqno (5.7)
$$
and for the remainder we have following estimate
$$
|R[\varphi]|\leq \|\varphi|L_2^{(2)}(0,1)\|\cdot
 0.01398. \eqno (5.8)
$$

For $N=3$ the optimal quadrature formula is
$$
\int\limits_0^1\varphi(x)dx=\frac{1}{30}\left[4\varphi(0)+11\varphi(1/3)+11\varphi(2/3)+
4\varphi(1)\right]+R[\varphi]\eqno (5.9)
$$
and for the remainder $R[\varphi]$ following estimate is valid
$$
|R[\varphi]|\leq \|\varphi|L_2^{(2)}(0,1)\|\cdot
 0.00586. \eqno (5.10)
$$

When $N=4$ the optimal quadrature formula is
$$
\int\limits_0^1\varphi(x)dx=\frac{1}{112}\left[11\varphi(0)+32\varphi(1/4)+26\varphi(1/2)+
32\varphi(3/4)+11\varphi(1)\right]+R[\varphi]\eqno (5.11)
$$
and for the remainder $R[\varphi]$ we get
$$
|R[\varphi]|\leq \|\varphi|L_2^{(2)}(0,1)\|\cdot 0.00305. \eqno
(5.12)
$$

For $N=5$ we obtain following optimal formula
$$
\int\limits_0^1\varphi(x)dx=\frac{1}{190}[15\varphi(0)+43\varphi(1/5)+37\varphi(2/5)+
$$
$$
+ 37\varphi(3/5)+43\varphi(4/5)+15\varphi(1)]+R[\varphi]\eqno
(5.13)
$$
with
$$
|R[\varphi]|\leq \|\varphi|L_2^{(2)}(0,1)\|\cdot 0.00188. \eqno
(5.14)
$$

Now we consider the quadrature formulas of the form (5.2).

From (5.2), (5.4), (5.5), (5.6) when $m=2$ and $\eta_0=0$ for
$N=2,3,4,5$ we get the optimal quadrature formulas (5.7)-(5.14),
respectively.

From (5.2), (5.4), (5.5), (5.6) for $m=2$ and $\eta_0=0.205$ when
$N=2,3,4,5$ we get following optimal quadrature formulas.

For $N = 2$ the optimal quadrature formula is
$$ \int\limits_0^1 {\varphi (x)dx = C_0 \varphi
(\eta _0 h) + C_1 } \varphi (1/2) + C_2 \varphi (1 - \eta _0 h) +
R[\varphi ],\eqno (5.15)
$$
where \\
\begin{tabular}{l}
$C_0=0.27075812274368231046,$\\
$C_1=0.45848375451263537906,$\\
$C_2=0.27075812274368231046$ \\
 \end{tabular} \\
and for the remainder $R[\varphi]$ following is valid
    $$
|R[\varphi ]| = |(\ell ,\varphi )| \le ||\varphi |L_2^{(2)}
(0,1)||_2  \cdot 0.00694814. \eqno (5.16)
$$

For $N = 3$ the optimal quadrature formula is
 $$
\int\limits_0^1 {\varphi (x)dx = C_0 \varphi (\eta _0 h) + C_1 }
\varphi (1/3) + C_2 \varphi (2/3) + C_3 \varphi (1 - \eta _0 h) +
R[\varphi ],\eqno (5.17)
$$
where \\
\begin{tabular}{l}
$C_0=0.17683465959328028293,$\\
$C_1=0.32316534040671971706,$ \\
$C_2=0.32316534040671971706,$\\
$C_3=0.17683465959328028293$\\
\end{tabular} \\
and  $$ |R[\varphi ]| = |(\ell ,\varphi )| \le ||\varphi
|L_2^{(2)} (0,1)||_2\cdot 0.00340515.  \eqno(5.18)
$$

For $N = 4$ the optimal formula is
    $$
\int\limits_0^1 {\varphi (x)dx = C_0 \varphi (\eta _0 h) + C_1 }
\varphi (1/4) + C_2 \varphi (1/2) + C_3 \varphi (3/4) + C_4
\varphi (1 - \eta _0 h) + R[\varphi ],\eqno (5.19)
$$
where \\ \begin{tabular}{l}
$C_0=0.13336566440349175557,$\\
$C_1=0.23884578079534432589,$\\
$C_2=0.25557710960232783705,$ \\
$C_3=0.23884578079534432589,$ \\
$C_4=0.13336566440349175557$\\
\end{tabular} \\
and
$$ |R[\varphi ]| = |(\ell ,\varphi )| \le ||\varphi
|L_2^{(2)} (0,1)||_2  \cdot 0.0020343. \eqno (5.20)
$$

For $N = 5$ the optimal quadrature formula is
    $$
\int\limits_0^1 {\varphi (x)dx = C_0 \varphi (\eta _0 h) + C_1 }
\varphi (1/5) + C_2 \varphi (2/5) + C_3 \varphi (3/5) +
$$
$$
+ C_4 \varphi (4/5) + C_5 \varphi (1 - \eta _0 h) + R[\varphi
],\eqno (5.21)
$$
\begin{tabular}{l}
$C_0=0.10653409090909090909,$\\
$C_1=0.19183238636363636363,$\\
$C_2=0.20163352272727272727,$\\
$C_3=0.20163352272727272727,$\\
$C_4=0.19183238636363636363,$\\
$C_5=0.10653409090909090909$\\
\end{tabular}\\
and for the remainder following estimate is valid
    $$
|R[\varphi ]| = |(\ell ,\varphi )| \le ||\varphi |L_2^{(2)}
(0,1)||_2  \cdot 0.0013408. \eqno (5.22)
$$

Thus, hence clear that the errors  (5.16), (5.18), (5.20), (5.22)
of the quadrature formulas (5.15), (5.17), (5.19), (5.21) are less
than the errors (5.8), (5.10), (5.12), (5.14) of the quadrature
formulas (5.7), (5.9), (5.11), (5.13), respectively.

\section{Acknowledgments}

The authors very thankful to professor G.V.Milovanovic for
discussion of the results and for some bibliographic references.

Kh.M.Shadimetov, A.R.Hayotov*\\
\emph{Institute of Mathematics and Information Technologies,\\
Uzbek Academy of Sciences,\\} e-mail*: hayotov@mail.ru,
abdullo\_hayotov@mail.ru

\end{document}